\documentclass[a4paper,12pt]{article}     

\usepackage{graphicx}		  
\usepackage{amsmath,amssymb}	

\newtheorem{lemma}{Lemma}[section]
\newtheorem{theorem}[lemma]{Theorem}
\newtheorem{proposition}[lemma]{Proposition}
\newtheorem{corollary}[lemma]{Corollary}
\newtheorem{assumption}{Assumption}

\def\authorfont{\footnotesize}

\def\ccode#1{\par
\vspace*{8pt}
{\authorfont{\leftskip18pt\rightskip\leftskip
\noindent #1\par}}\par}

\newenvironment{Proof}{
\hspace*{-9mm}
{ \it Proof.}}
{\hfill {$\square$}\vspace{1.5em}}

\begin{document}

\begin{center}{
{\Large 
 Properties of minimal charts and
 their applications XI: 
 no minimal charts with exactly seven white vertices}
\vspace{10pt}
\\ 
Teruo NAGASE and Akiko SHIMA\footnote{The second author is supported by JSPS KAKENHI Grant Number 21K03255.}
}
\end{center}

\begin{abstract}
Charts are oriented labeled graphs in a disk.
Any simple surface braid (2-dimensional braid) can be described by using a chart.
Also, a chart represents an oriented closed surface
embedded in 4-space.
In this paper, we investigate embedded surfaces in 4-space
by using charts.
In this paper, we shall show that 
there is no minimal chart with exactly seven white vertices.
\end{abstract}

%
%
%
%

\ccode{2020 Mathematics Subject Classification. Primary 57K45,05C10; Secondary 57M15.}
\ccode{ {\it Key Words and Phrases}. surface link, chart, C-move, white vertex. }


\setcounter{section}{0}
\section{Introduction}


Charts are oriented labeled graphs in a disk (see  \cite{KnottedSurfaces},\cite{BraidBook}, and see Section~\ref{s:Prel}  for the precise definition of charts).
Let $D_1^2, D_2^2$ be 2-dimensional disks.
Any simple surface braid (2-dimensional braid) can be described 
by using a chart,
here a simple surface braid is a properly embedded surface
$S$ in the 4-dimensional disk $D_1^2\times D_2^2$ such that
a natural map $\pi:S\subset D_1^2\times D_2^2\to D_2^2$ is 
a simple branched covering map of $D_2^2$ and
the boundary $\partial S$ is a trivial closed braid in
the solid torus $D_1^2\times \partial D_2^2$
(see \cite{BraidThree}, \cite[Chapter 14 and Chapter 18]{BraidBook}).
Also, from a chart, 
we can construct a simple closed surface braid in 4-space ${\Bbb R}^4$. This surface is an oriented closed surface 
embedded in ${\Bbb R}^4$.
On the other hand, any oriented embedded closed surface 
 in ${\Bbb R}^4$ is ambient isotopic to a simple
closed surface braid
 (see \cite{BraidThree},\cite[Chapter 23]{BraidBook}). 
A C-move 
is a local modification between two charts
in a disk (see Section~\ref{s:Prel} for C-moves).
A C-move between two charts induces 
an ambient isotopy between oriented closed surfaces 
corresponding to the two charts.
In this paper, we investigate oriented closed surfaces in 4-space
by using charts.

We will work in the PL category or smooth category. All submanifolds are assumed to be locally flat.
In \cite{ONS},
we showed that there is no minimal chart with exactly five vertices
 (see Section~\ref{s:Prel} for the precise definition of minimal charts). 
Hasegawa proved that there exists a minimal chart with exactly
six white vertices \cite{H1}. 
This chart represents a 2-twist spun trefoil.
In \cite{INS} and \cite{NST},
we investigated minimal charts with exactly four white vertices.
In this paper, 
we investigate properties of minimal charts 
which show that
there is no minimal chart with exactly seven white vertices
(see \cite{ChartApp1},\cite{ChartAppII},\cite{ChartAppIII},\cite{ChartAppIV},\cite{ChartAppV},\cite{ChartAppVI},\cite{ChartAppVII},\cite{ChartAppVIII},\cite{ChartAppIX},\cite{ChartAppX}).

Let $\Gamma$ be a chart.
For each label $m$, we denote by $\Gamma_m$
the union of all the edges of label $m$.

Now we define a type of a chart:
Let $\Gamma$ be a chart with at least one white vertex, 
and $n_1,n_2,\dots,n_k$ integers.
The chart $\Gamma$ is of {\it type $(n_1,n_2,\dots,n_k)$} if there exists a label $m$ of $\Gamma$ satisfying the following three conditions:
\begin{enumerate}
\item[(i)] For each $i=1,2,\dots, k$, 
the chart $\Gamma$ contains exactly $n_{i}$ white vertices in $\Gamma_{m+i-1}\cap \Gamma_{m+i}$.
\item[(ii)] If $i<0$ or $i>k$, then $\Gamma_{m+i}$ does not contain any white vertices.
\item[(iii)] Both of the two subgraphs $\Gamma_m$ and $\Gamma_{m+k}$ contain at least one white vertex.
\end{enumerate}
If we want to emphasize the label $m$,
then we say that $\Gamma$ is of {\it type $(m;n_1,n_2,\dots,n_k)$}. 
Note that $n_1\ge1$ and $n_k\ge1$ by Condition~(iii).

We proved in \cite[Theorem 1.1]{ChartAppII} that
if there exists a minimal $n$-chart $\Gamma$ with exactly seven white vertices,
then $\Gamma$ is a chart of 
type $(7),(5,2),(4,3),(3,2,2)$ or $(2,3,2)$ 
(if necessary we change the label
$i$ by $n-i$ for all label $i$).
In \cite{ChartAppV},
we showed that
there is no minimal chart of type $(3,2,2)$.
In \cite{ChartAppVI} and \cite{ChartAppVII},
there is no minimal chart of type $(2,3,2)$.
In \cite{ChartAppVIII},
there is no minimal chart of type $(7)$.
In \cite{ChartAppIX},
there is no minimal chart of type $(4,3)$.
In \cite{ChartAppX}, 
we investigate a minimal chart of type $(5,2)$.

In this paper we shall show the following:

\begin{theorem}
\label{MainTheorem} 
There is no minimal chart of type $(5,2)$.
\end{theorem}

From the above theorem, we have the following:

\begin{theorem}
\label{MainTheorem2} 
There is no minimal chart with exactly seven white vertices.
\end{theorem}

The paper is organized as follows.
In Section~\ref{s:Prel},
we define charts and minimal charts.
In Section~\ref{s:kAngledDisk},
 we review lemmas of a $2$-angled disk and
a $3$-angled disk of $\Gamma_m$
for a minimal chart $\Gamma$ and a label $m$, 
where
a $k$-angled disk is a disk 
whose boundary contains exactly $k$ white vertices
and consists of edges of label $m$.
In Section~\ref{s:ThetaCurve},
 we shall show that if $\Gamma$ is a minimal chart 
of type $(m; 5, 2)$,
then the graph $\Gamma_m$ does not 
contain a $\theta$-curve as shown in Fig.~\ref{fig09}(a).
In Section~\ref{s:DiskLemma},
we review a useful lemma called New Disk Lemma(Lemma~\ref{NewDiskLemma}), and we shall extend this lemma.
In Section~\ref{s:IOC},
we review IO-Calculation(a property of numbers of 
inward arcs of label $k$ 
and outward arcs of label $k$ in a closed domain $F$
with $\partial F\subset\Gamma_{k-1}\cup\Gamma_k\cup\Gamma_{k+1}$
for some label $k$).
In Section~\ref{s:Lens},
we review a useful lemma for a disk called a lens.
In Section~\ref{s:Oval},
 we shall show that if $\Gamma$ is a minimal chart 
of type $(m; 5, 2)$,
then the graph $\Gamma_m$ does not contain
an oval as shown 
in Fig.~\ref{fig09}(b).
In Section~\ref{s:TypeH},
 we shall show that if $\Gamma$ is a minimal chart 
of type $(m; 5, 2)$,
then the graph $\Gamma_m$ does not contain
the graph as shown 
in Fig.~\ref{fig10}(h).
In Section~\ref{s:TypeG},
 we shall show that if $\Gamma$ is a minimal chart 
of type $(m; 5, 2)$,
then the graph $\Gamma_m$ does not contain
the graph as shown 
in Fig.~\ref{fig10}(g).
Moreover, we shall prove Theorem~\ref{MainTheorem}.


\section{Preliminaries}
\label{s:Prel}

In this section, 
we introduce 
the definition of charts and its related words.

\begin{figure}[htb]
\begin{center}
\includegraphics{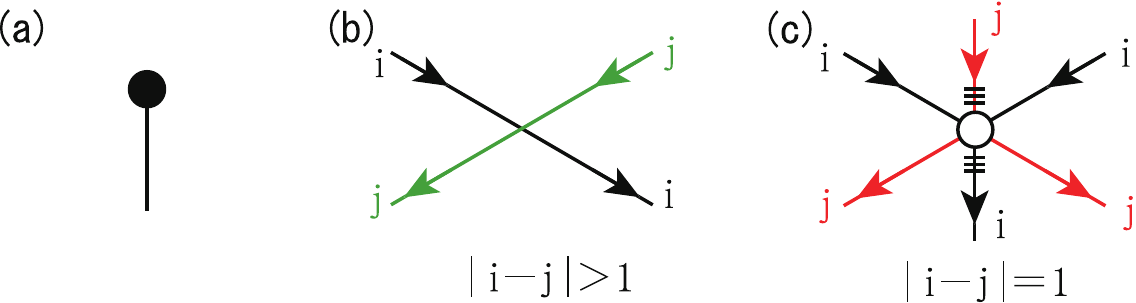}
\end{center}
\caption{ \label{fig01} (a) A black vertex. (b) A crossing. (c) A white vertex. 
Each arc with three transversal short arcs is a middle arc at the white vertex. }
\end{figure}

Let $n$ be a positive integer.
An $n$-{\it chart}  
(a braid chart of degree $n$ \cite{KnottedSurfaces}
or a surface braid chart of degree $n$ \cite{BraidBook}) 
is 
an oriented labeled graph in the interior of a disk,
which may be empty 
or
have closed edges without vertices
satisfying the following four conditions
(see Fig.~\ref{fig01}):
\begin{enumerate}
\item[(i)] 
Every vertex has degree $1$, $4$, or $6$.
\item[(ii)] 
The labels of edges are 
in $\{1,2,\dots,n-1\}$.
\item[(iii)]
In a small neighborhood of
each vertex of degree $6$,
there are six short arcs,
three consecutive arcs are
oriented inward 
and
the other three are outward,
and
these six are labeled $i$ and $i+1$
alternately for some $i$,
where the orientation and label of
each arc are inherited from
the edge containing the arc.
\item[(iv)]
For each vertex of degree $4$,
diagonal edges have the same label
and
are oriented coherently,
and the labels $i$ and $j$ of
the diagonals satisfy $|i-j|>1$.
\end{enumerate}
We call a vertex of degree $1$ a {\it black vertex},
a vertex of degree $4$ a {\it crossing}, and 
a vertex of degree $6$ a {\it white vertex}
respectively.

Among six short arcs
in a small neighborhood of
a white vertex,
a central arc of each three consecutive arcs
oriented inward (resp. outward) 
is called a   
{\it middle arc} at the white vertex
(see Fig.~\ref{fig01}(c)).
For each white vertex $v$, 
there are two middle arcs at $v$ 
in a small neighborhood of $v$.
An edge is said to be {\it middle at} a white vertex $v$ if it contains a middle arc at $v$.

Let $e$ be an edge connecting $v_1$ and $v_2$.
If $e$ is oriented from $v_1$ to $v_2$,
then we say that 
$e$ is oriented {\it outward at $v_1$}
and {\it inward at $v_2$}.


\begin{figure}
\begin{center}
\includegraphics{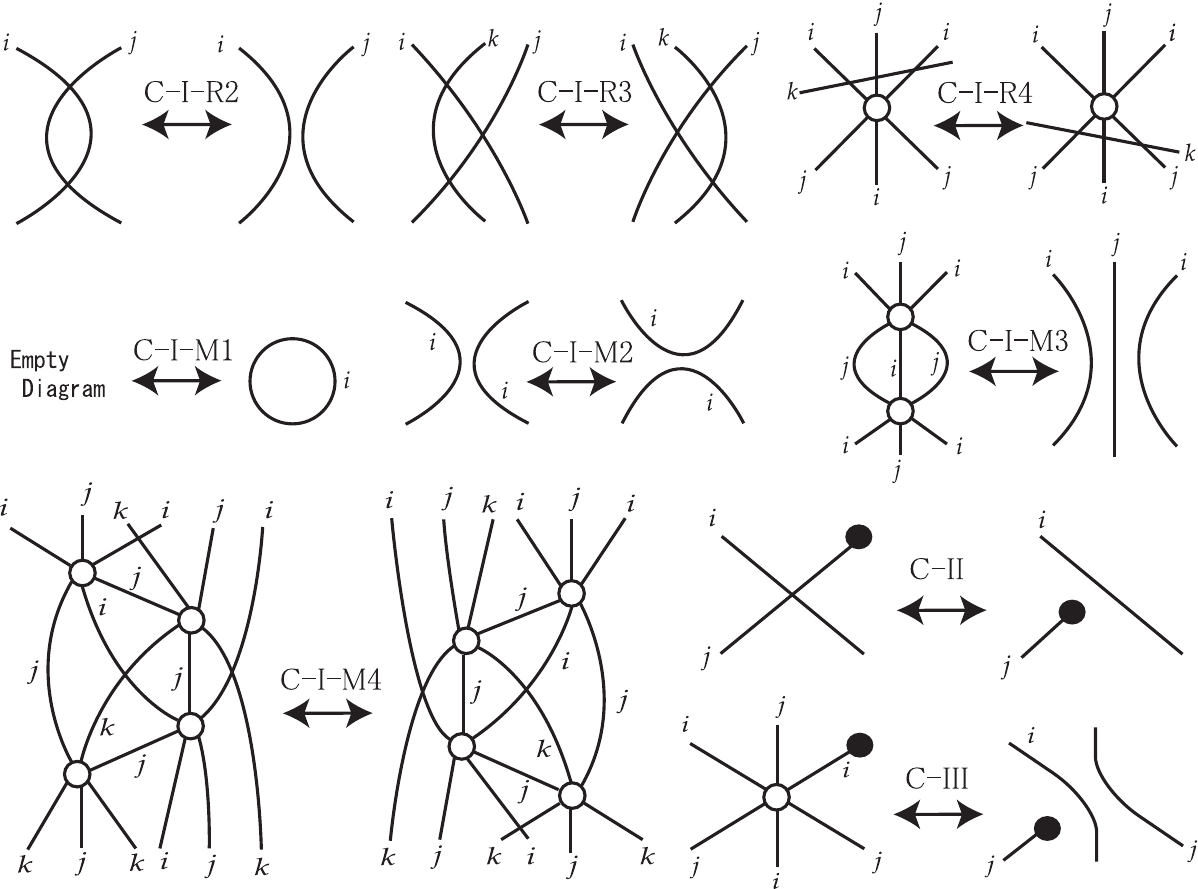}
\end{center}
\caption{ \label{fig02} For the C-III move, 
the edge with the black vertex is not middle at
a white vertex in the left figure. }
\end{figure}

Now {\it C-moves} are local modifications 
of charts as shown in Fig.~\ref{fig02}
(cf. \cite{KnottedSurfaces}, 
\cite{BraidBook} and \cite{Tanaka}).
Two charts are said to be {\it C-move equivalent}  if there exists
a finite sequence of C-moves 
which modifies one of the two charts 
to the other.

An edge in a chart is called 
a {\it free edge}
if it has
two black vertices.

For each chart $\Gamma$,
let $w(\Gamma)$ and $f(\Gamma)$ be the number of white vertices, and the number of free edges respectively.
The pair $(w(\Gamma), -f(\Gamma))$ is called a {\it complexity} of the chart (see \cite{BraidThree}).
A chart $\Gamma$ is called a {\it minimal chart} if its complexity is minimal among the charts C-move equivalent to the chart $\Gamma$ with respect to the lexicographic order of pairs of integers.

We showed the difference of a chart in a disk and in a 2-sphere (see \cite[Lemma 2.1]{ChartApp1}).
This lemma follows from that there exists a natural one-to-one correspondence between $\{$charts in $S^2\}/$C-moves and $\{$charts in $D^2\}/$C-moves, conjugations
(\cite[Chapter 23 and Chapter 25]{BraidBook}).
To make the argument simple, we assume that 
the charts lie on the 2-sphere instead of the disk.
\begin{assumption}
In this paper,
all charts are contained in the $2$-sphere $S^2$.
\end{assumption}
We have the special point in the 2-sphere $S^2$, called the point at infinity,
 denoted by $\infty$.
In this paper, all charts are contained in a disk such that the disk 
does not contain the point at infinity $\infty$.

Let $\Gamma$ be a chart,
and $m$ a label of $\Gamma$. 
A {\it hoop} is a closed edge of $\Gamma$ without vertices 
(hence without crossings, neither).
A {\it ring} is a simple closed curve in $\Gamma_m$ containing at least one crossing but not containing any white vertices.
A hoop is said to be {\it simple} 
if one of the two complementary domains
of the hoop
does not contain any white vertices.

An edge in a chart is called 
a {\it terminal edge}
if it has
a white vertex and a black vertex.

We can assume that
all minimal charts $\Gamma$
satisfy the following four conditions 
(see \cite{ChartApp1},\cite{ChartAppII},\cite{ChartAppIII},\cite{StI}):

\begin{assumption}
\label{AssumeTerminal}
If an edge of $\Gamma$
contains a black vertex,
then the edge is a free edge 
or a terminal edge.
Moreover 
any terminal edge contains a middle arc.
\end{assumption}

\begin{assumption}
\label{NoSimpleHoop}
All free edges and simple hoops in $\Gamma$ 
are moved into a small neighborhood $U_\infty$ 
of the point at infinity $\infty$. 
Hence
we assume that 
$\Gamma$ does not contain free edges
nor simple hoops, 
otherwise mentioned. 
\end{assumption}

\begin{assumption}
\label{Ring}
Each complementary domain of
any ring and hoop must contain 
at least one white vertex. 
\end{assumption}

\begin{assumption}
\label{Infinity}
The point at infinity $\infty$ is moved in any complementary domain of $\Gamma$.
\end{assumption}

In this paper
for a subset $X$ in a space
we denote 
the interior of $X$,
the boundary of $X$ and
the closure of $X$
by Int$X$, $\partial X$
and $Cl(X)$
respectively.

Let $\alpha$ be a simple arc or an edge,
and $p,q$ the endpoints of $\alpha$.
We denote 
$\partial \alpha=\{p,q\}$
and ${\rm Int}\alpha=\alpha-\{p,q\}$.


\section{$k$-angled disks}
\label{s:kAngledDisk}

In this section, we review lemmas for a disk called a $k$-angled disk.

Let $\Gamma$ be a chart, $m$ a label of $\Gamma$, 
$D$ a disk with $\partial D\subset \Gamma_m$, 
and $k$ a positive integer.
If $\partial D$ contains exactly
$k$ white vertices, 
then $D$ is called 
{\it a $k$-angled disk of $\Gamma_m$}. 
Note that 
the boundary $\partial D$ may contain crossings.

Let $\Gamma$ be a chart, and
$m$ a label of $\Gamma$.
An edge of label $m$ is called a {\it feeler} of a $k$-angled disk $D$ of $\Gamma_m$
if the edge intersects $N-\partial D$
where $N$ is a regular neighborhood of $\partial D$ in $D$.

Let $\Gamma$ be a chart. 
Suppose that an object consists of 
some edges of $\Gamma$, arcs in edges of 
$\Gamma$ and arcs around white vertices.
Then the object is called a {\it pseudo chart}.

Let $X$ be a set in a chart $\Gamma$.
Let
 $$w(X)=\text{the number of white vertices in $X$.}$$

\begin{figure}
\centerline{\includegraphics{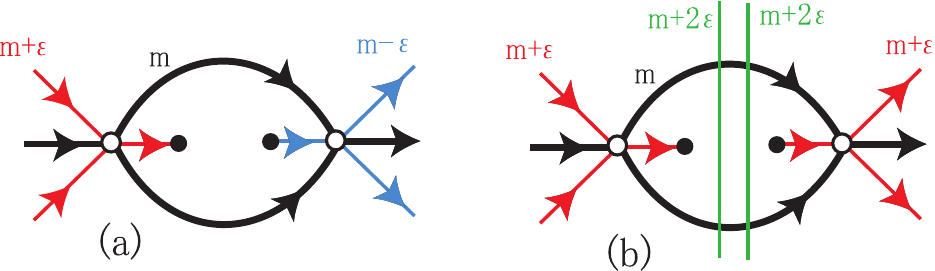}}
\caption{\label{fig03}
$m$ is a label,
and $\varepsilon\in\{+1,-1\}$.}
\end{figure}

\begin{lemma}
\label{Theorem2AngledDisk}
{\rm (\cite[Corollary 5.8]{ChartAppII})}
Let $\Gamma$ be a minimal chart.
Let $D$ be a $2$-angled disk of $\Gamma_m$ with at most one feeler.
If $w(\Gamma\cap{\rm Int}D)=0$,
then a regular neighborhood of $D$ contains one of two pseudo charts as shown in Fig.~\ref{fig03}.
\end{lemma}

Let $\Gamma$ be a chart, 
$D$ a disk, and 
$G$ a pseudo chart with $G \subset D$.
Let $r:D\to D$ be a reflection of $D$, and $G^*$ the pseudo chart obtained from $G$ by changing the orientations of all of the edges.
Then the set $\{G,G^*, r(G), r(G^*)\}$ 
is called the {\it RO-family of the pseudo chart $G$}.

Let $\Gamma$ be a chart,
and $D$ a $k$-angled disk of $\Gamma_m$.
If any feeler of $D$ of label $m$ is a terminal edge,
then $D$ is called a {\it special} $k$-angled disk.

\begin{figure}[htb]
\centerline{\includegraphics{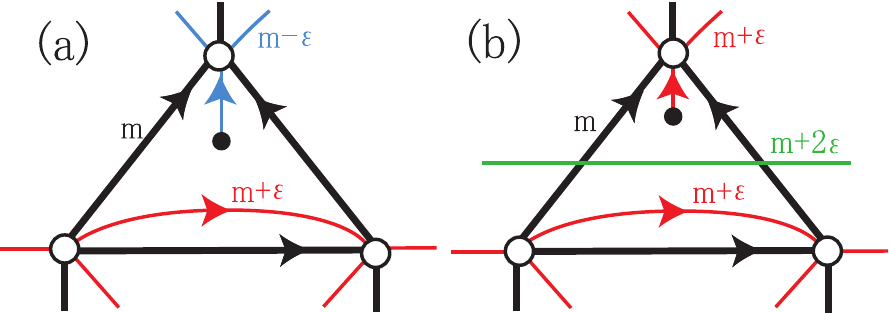}}
\caption{\label{fig04}
The 3-angled disks  have no feelers,
$m$ is a label, $\varepsilon\in\{+1,-1\}$.}
\end{figure}

\begin{lemma}
{\rm (\cite[Lemma 4.2(a)]{ChartAppIX})}
\label{Theorem3AngledDisk}
Let $\Gamma$ be a minimal chart, and $m$ a label of $\Gamma$.
Let $D$ be a special $3$-angled disk of $\Gamma_m$
with at most two feelers.
If $w(\Gamma\cap {\rm Int}D)=0$,
then a regular neighborhood of $D$ contains one of the RO-families of the two pseudo charts as shown in 
Fig.~\ref{fig04}.

\end{lemma}

\begin{figure}[htb]
\centerline{\includegraphics{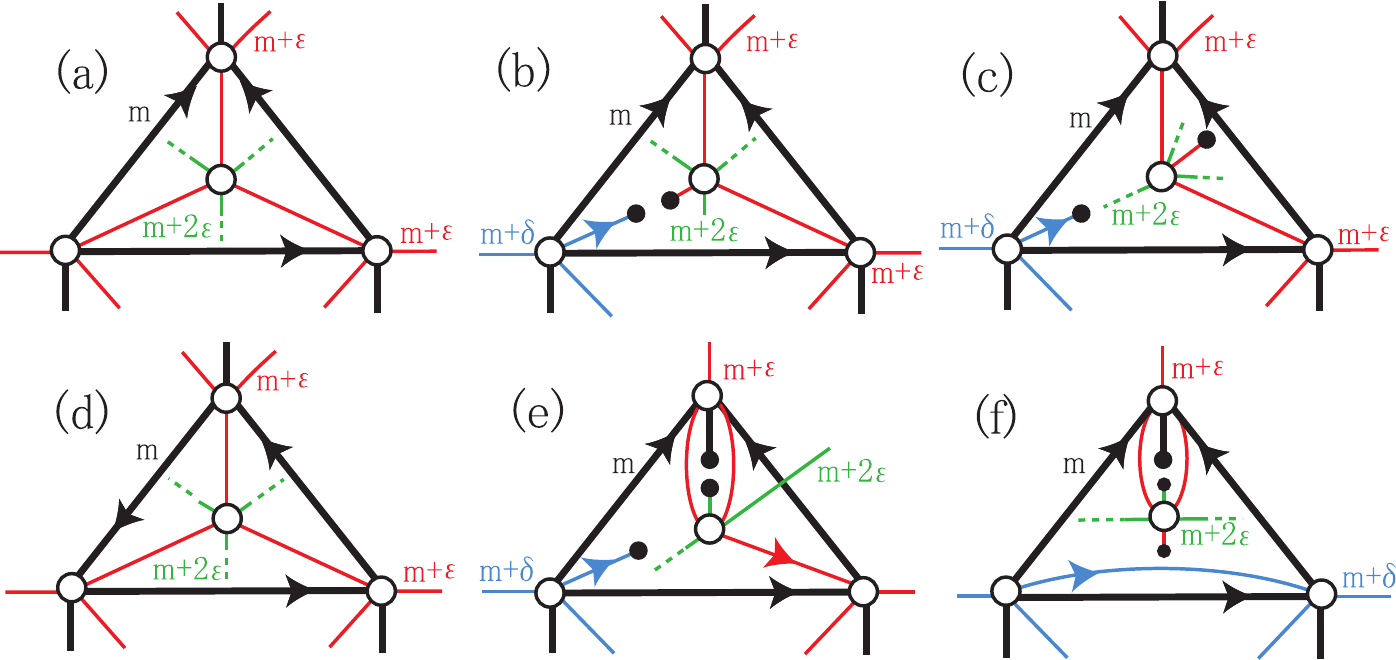}}
\caption{\label{fig05}
(a),(b),(c),(d) 3-angled disks without feelers.
(e),(f) 3-angled disks with one feeler.}
\end{figure}

\begin{lemma}
{\rm (\cite[Lemma 4.2(b)]{ChartAppIX})}
\label{Theorem3AngledDiskChart9}
Let $\Gamma$ be a minimal chart, and $m$ a label of $\Gamma$.
Let $D$ be a special $3$-angled disk of $\Gamma_m$
with at most two feelers.
If $w(\Gamma\cap {\rm Int}D)=w(\Gamma_{m+\varepsilon}\cap {\rm Int}D)=1$ for some $\varepsilon\in\{+1,-1\}$,
then a regular neighborhood of $D$ contains one of the RO-families of the six pseudo charts as shown in 
Fig.~\ref{fig05}.
\end{lemma}

\begin{figure}[htb]
\begin{center}
\centerline{\includegraphics{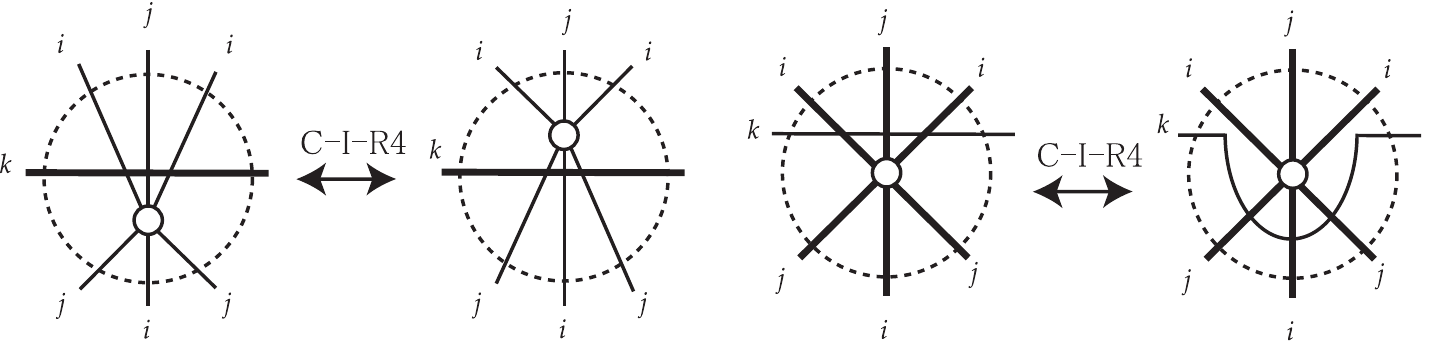}}
\caption{\label{fig06}
C-moves keeping thicken figures fixed.}
\end{center}
\end{figure}

Let $\Gamma $ and $\Gamma^\prime $ be C-move equivalent charts. 
Suppose that a pseudo chart $X$ of $\Gamma$ is also a pseudo chart of $\Gamma^\prime$. 
Then we say that 
$\Gamma$ is modified to $\Gamma^\prime$ by {\it C-moves keeping $X$ fixed}.
In Fig.~\ref{fig06},
we give examples of C-moves keeping pseudo charts  fixed.

Let $\Gamma$ be a chart, and $X$ a subset of $\Gamma$. Let 
$$c(X)=\text{the number of crossings in $X$}.$$

Let $D$ be a $k$-angled disk of $\Gamma_m$ for a minimal chart $\Gamma$.
The pair of integers 
$(w(\Gamma\cap {\rm Int}D),c(\partial D))$
is called the {\it local complexity 
with respect to $D$},
denoted by $\ell c(D;\Gamma)$.
Let
${\Bbb S}$ be the set of all minimal charts each of which can be moved from $\Gamma$ by C-moves in a regular neighborhood of $D$ keeping $\partial D$ fixed.
The chart $\Gamma$ is said to be 
{\it locally minimal
with respect to $D$}
if its local complexity
with respect to $D$
is minimal
among the charts in ${\Bbb S}$ with respect to 
the lexicographic order.

\begin{figure}[htb]
\begin{center}
\centerline{\includegraphics{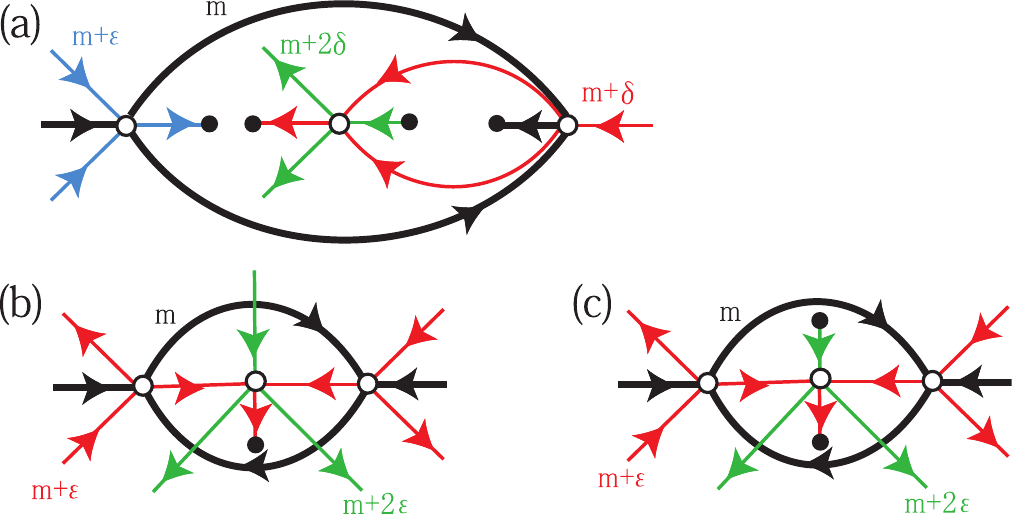}}
\caption{\label{fig07}
The 2-angled disk (a) has one feeler, the others do not have any feelers.} 
\end{center}
\end{figure}

\begin{lemma}
\label{Theorem2AngledDisk2}
{\rm (\cite[Theorem 1.1]{ChartAppIII})}
Let $\Gamma$ be a minimal chart.
Let $D$ be a $2$-angled disk of $\Gamma_m$ with at most one feeler
such that $\Gamma$ is locally minimal
with respect to $D$.
If $w(\Gamma\cap{\rm Int}D)\leqq1$,
then a regular neighborhood of $D$ contains an element in the RO-families of the five pseudo charts as shown in 
Fig.~\ref{fig03} and Fig.~\ref{fig07}.
\end{lemma}

\section{Case of the $\theta$-curve}
\label{s:ThetaCurve}

In this section, we shall show that if 
$\Gamma$ is a minimal chart of type $(m;5,2)$,
then the graph $\Gamma_m$ does not contain a $\theta$-curve.

\begin{figure}[htb]
\centerline{\includegraphics{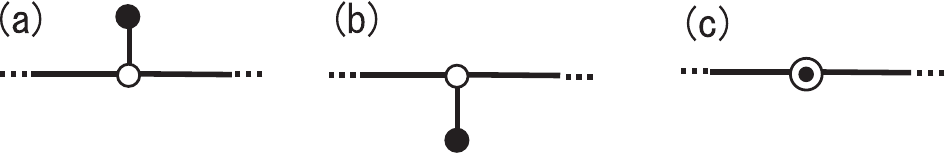}}
\caption{\label{fig08}
(a),(b) White vertices in terminal edges.
(c) BW-vertex.}
\end{figure}

 In our argument  we often construct a chart $\Gamma$. 
On the construction of a chart $\Gamma$, for a white vertex $w\in\Gamma_m$ for some label $m$,  
among the three edges of $\Gamma_m$ 
containing $w$, 
if one of the three edges is a terminal edge 
(see Fig.~\ref{fig08}(a) and (b)), 
then we remove the terminal edge and
put a black dot at the center of the white vertex  as shown in Fig.~\ref{fig08}(c).
Namely
Fig.~\ref{fig08}(c) means 
Fig.~\ref{fig08}(a) or 
Fig.~\ref{fig08}(b).
We call the vertex in Fig.~\ref{fig08}(c) 
a {\it BW-vertex} with respect to $\Gamma_m$.

\begin{figure}
\centerline{\includegraphics{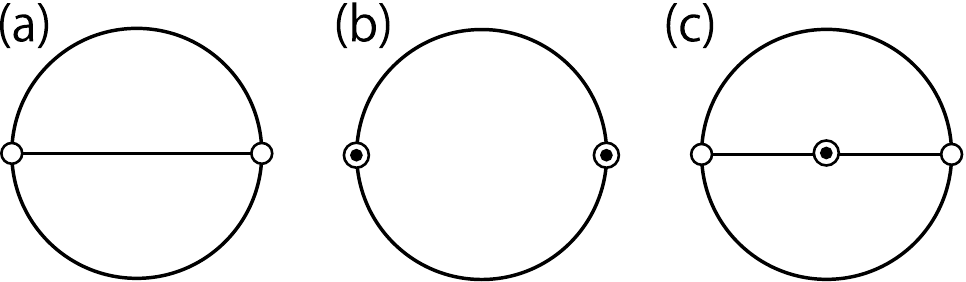}}
\caption{\label{fig09}
(a) A $\theta$-curve. (b) An oval. 
(c) A skew $\theta$-curve.}
\end{figure}

The three graphs in Fig.~\ref{fig09}
are examples of graphs in $\Gamma_m$ for a chart $\Gamma$
and a label $m$.
We call 
a {\it $\theta$-curve},
an {\it oval},
a {\it skew $\theta$-curve} the three graphs as shown
in Fig.~\ref{fig09}(a),(b),(c)
respectively.

Let $\Gamma$ be a chart,
and $m$ a label of $\Gamma$. 
A {\it loop} is a simple closed curve in $\Gamma_m$ with exactly one white vertex
(possibly with crossings).

\begin{figure}
\centerline{\includegraphics{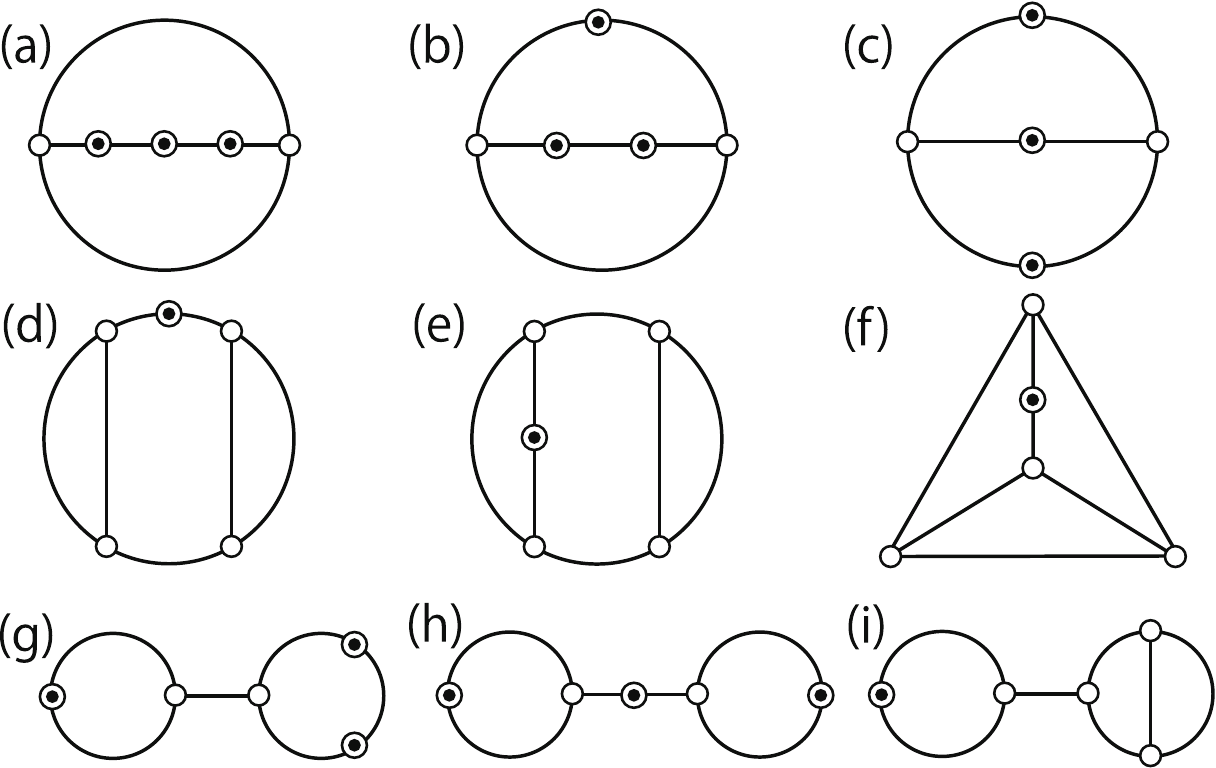}}
\caption{\label{fig10}
(a),(b),(c),(g),(h) Graphs with three black vertices.
(d),(e),(f),(i) Graphs with one black vertex.}
\end{figure}

\begin{lemma}
\label{GammaMFiveWhite}
{\rm (\cite[Lemma 3.5]{ChartAppVI})}
Let $\Gamma$ be a minimal chart,
and $m$ a label of $\Gamma$.
If $w(\Gamma_m)=5$ and if $\Gamma_m$ has no loop,
then the graph $\Gamma_m$ contains one of the following graphs:
\begin{enumerate}
\item[{\rm (a)}] one of the nine graphs as shown in 
Fig.~\ref{fig10},
\item[{\rm (b)}] the union of a $\theta$-curve and a skew $\theta$-curve,
\item[{\rm (c)}] the union of an oval and a skew $\theta$-curve.
\end{enumerate}
\end{lemma}

\begin{lemma}$(${\rm \cite[Theorem 1.1]{ChartAppIV}}$)$
\label{LemmaNoLoop}
There is no loop in any minimal chart with exactly seven white vertices.
\end{lemma}

Let $\Gamma$ be a chart, 
and $m$ a label of $\Gamma$.
Let $L$ be the closure of a connected component 
of the set obtained by taking out 
all the white vertices from $\Gamma_m$.
If $L$ contains at least one white vertex
but does not contain any black vertex,
then $L$ is called an {\it internal edge of label $m$}.
Note that an internal edge may contain a crossing of $\Gamma$.

\begin{lemma}
\label{NoThetaCurve}
Let $\Gamma$ be a minimal chart of type $(m;5,2)$.
Then $\Gamma_m$ does not contain a $\theta$-curve.
\end{lemma}

\begin{figure}[htb]
\centerline{\includegraphics{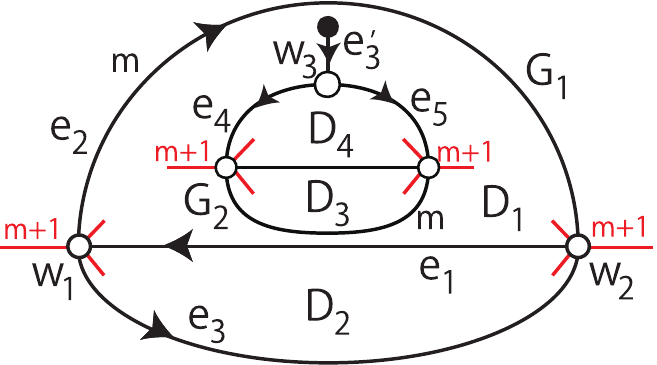}}
\caption{\label{fig11}
A $\theta$-curve and a skew $\theta$-curve.}
\end{figure}

\begin{Proof}
By Lemma~\ref{LemmaNoLoop}, the chart $\Gamma$ has no loop.
Hence $\Gamma_m$ has no loop.

Suppose that $\Gamma_m$ contains a $\theta$-curve, say $G_1$.
Then by Lemma~\ref{GammaMFiveWhite},
the graph $\Gamma_m$ contains a skew $\theta$-curve, say $G_2$.
Let $e_1,e_2,e_3$ be the three internal edges 
 of label $m$ in $G_1$,
and $w_1,w_2$ the white vertices in $G_1$.
Without loss of generality,
we can assume that 
\begin{enumerate}
\item[(1)] $e_1$ is oriented inward at $w_1$ and middle at $w_1$.
\end{enumerate}
Then 
\begin{enumerate}
\item[(2)] $e_2,e_3$ are oriented outward at $w_1$,
\item[(3)] $e_1$ is middle at $w_2$.
\end{enumerate}

Let $D_1,D_2$ be the special 2-angled disks of $\Gamma_m$
with $\partial D_1=e_1\cup e_2$ and
$\partial D_2=e_1\cup e_3$.
Then by (1) and (2),
both of $\partial D_1$ and $\partial D_2$
are oriented clockwise or anticlockwise.
Moreover, by (1) and (3), the edge $e_1$ is middle at 
both white vertices $w_1$ and
$w_2$.
Thus by Lemma~\ref{Theorem2AngledDisk2}, we have
\begin{enumerate}
\item[(4)] $w(\Gamma\cap {\rm Int}D_1)\geqq2$ and
$w(\Gamma\cap {\rm Int}D_2)\geqq2$.
\end{enumerate}

Since $G_1\cap G_2=\emptyset$,
we have $G_2\subset S^2-G_1$.
Hence $G_2\subset {\rm Int}D_1$ or $G_2\subset {\rm Int}D_2$
or $G_2\subset S^2-(D_1\cup D_2)$.

We shall show that $G_2\subset {\rm Int}D_1$ or $G_2\subset {\rm Int}D_2$.
If $G_2\subset S^2-(D_1\cup D_2)$,
then \vspace{2mm}

$\begin{array}{rcl}
7=w(\Gamma) & \geqq & w(G_1)+w(G_2)+w(\Gamma\cap{\rm Int}D_1)+
w(\Gamma\cap{\rm Int}D_2)\vspace{2mm}\\
& \geqq & 2+3+2+2=9.\vspace{2mm}
\end{array}
$\\
This is a contradiction.
Hence $G_2\subset {\rm Int}D_1$ or $G_2\subset {\rm Int}D_2$.

Without loss of generality we can assume that 
$G_2\subset {\rm Int}D_1$.
Then the graph $G_2$ separates the disk $D_1$ into three regions.
Two of the three regions are disks, say $D_3,D_4$.
Note that $D_3,D_4$ are a 2-angled disk or a 3-angled disk.

We shall show that neither $D_3$ nor $D_4$ has a feeler.
If one of $D_3,D_4$ has a feeler,
then the disk is a 3-angled disk with exactly one feeler.
Thus by Lemma~\ref{Theorem3AngledDisk},
the disk contains at least one white vertex in its interior.
Hence $w(\Gamma\cap {\rm Int}D_3)\geqq1$ or
$w(\Gamma\cap {\rm Int}D_4)\geqq1$.
Thus $w(\Gamma\cap {\rm Int}D_1)\geqq 4$.
Hence by (4) \vspace{2mm}

$\begin{array}{rcl}
7=w(\Gamma) & \geqq & w(G_1)+w(\Gamma\cap{\rm Int}D_1)+
w(\Gamma\cap{\rm Int}D_2)\vspace{2mm}\\
& \geqq & 2+4+2=8.\vspace{2mm}
\end{array}
$\\
This is a contradiction.
Thus neither $D_3$ nor $D_4$ has a feeler
(see Fig.~\ref{fig11}).

Without loss of generality, we can assume that
$D_3$ is a 2-angled disk and
$D_4$ is a 3-angled disk.
Let $e_3'$ be the terminal edge of label $m$ in $G_2$,
and $w_3,w_4,w_5$ the white vertices in $G_2$
with $w_3\in e_3'$.
Let $e_4,e_5$ be internal edges of label $m$ in $G_2$
with $\partial e_4=\{w_3,w_4\}$ and
$\partial e_5=\{w_3,w_5\}$.

If necessary we change the orientation of all edges,
we can assume that
the terminal edge $e_3'$ is oriented inward at $w_3$.
Then by Assumption~\ref{AssumeTerminal},
both of $e_4,e_5$ are oriented outward at $w_3$.
Thus $e_4$ is oriented inward at $w_4$ and
$e_5$ is oriented inward at $w_5$.
Hence by Lemma~\ref{Theorem2AngledDisk},
we have $w(\Gamma\cap{\rm Int}D_3)\geqq1$.
However we can show $w(\Gamma)\geqq8$ 
by the similar way as above.
This contradicts the fact that $w(\Gamma)=7$.
Therefore,
the graph $\Gamma_m$ does not contain a $\theta$-curve. 
\end{Proof}

\section{Disk Lemma}
\label{s:DiskLemma}

In this section, we review a useful lemma called New Disk Lemma.
Moreover, we shall extend this lemma in this section.

Let $\Gamma$ be a chart, and $D$ a disk.
Let $\alpha$ be a simple arc in $\partial D$,
and $\gamma$ a simple arc in an internal edge of label $k$.
We call the simple arc $\gamma$
a {\it {$(D,\alpha)$-arc}} of label $k$
provided that 
$\partial \gamma \subset $Int$\alpha$
and
Int$\gamma\subset $Int$D$. 
If there is no $(D,\alpha)$-arc in $\Gamma$,
then the chart $\Gamma$ is said to be
$(D,\alpha)$-{\it arc free}.

\begin{lemma}
$($New Disk Lemma$)$
{\em (\cite[Lemma 7.1(a)]{Chart33},
 cf. \cite[Lemma 3.2]{ChartApp1})} 
\label{NewDiskLemma}
Let $\Gamma$ be a chart and
$D$ a disk 
whose interior does not contain 
a white vertex nor a black vertex of $\Gamma$.
Let $\alpha$ be a simple arc in $\partial D$ 
such that ${\rm Int}\alpha$ does not contain 
a white vertex nor a black vertex of $\Gamma$.
Let $V$ be a regular neighborhood of $\alpha$. 
If the arc $\alpha$ is contained in 
an internal edge of some label $k$ of $\Gamma$,
then by applying C-I-M2 moves, C-I-R2 moves, 
and C-I-R3 moves in $V$, 
there exists 
a $(D,\alpha)$-arc free chart $\Gamma'$ 
obtained from the chart $\Gamma$ 
keeping $\alpha$ fixed 
$($cf. Fig.~\ref{fig13}$)$.
\end{lemma}

Let $D$ be a disk, $\alpha$ and $\beta$ two simple arcs
with $\partial D=\alpha \cup \beta$, and
$\alpha\cap \beta=\partial \alpha=\partial \beta$.
The pair $(\alpha,\beta)$
is called a {\it boundary arc pair} of the disk $D$.

\begin{figure}[htb]
\centerline{\includegraphics{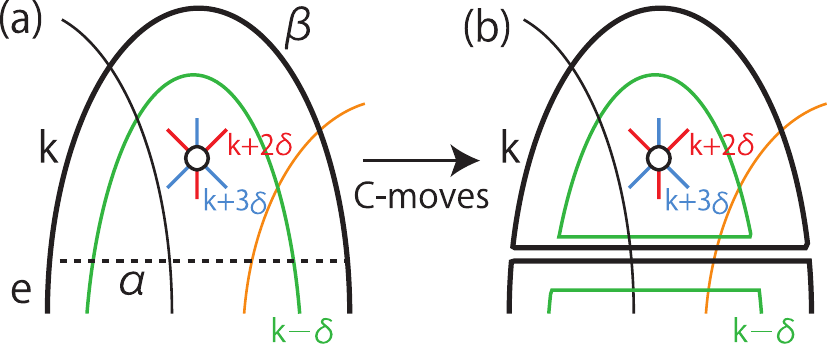}}
\caption{\label{fig12}
The edge $e$ can be moved the set $(e-\beta)\cup\alpha$ by C-moves. }
\end{figure}

\begin{figure}[htb]
\centerline{\includegraphics{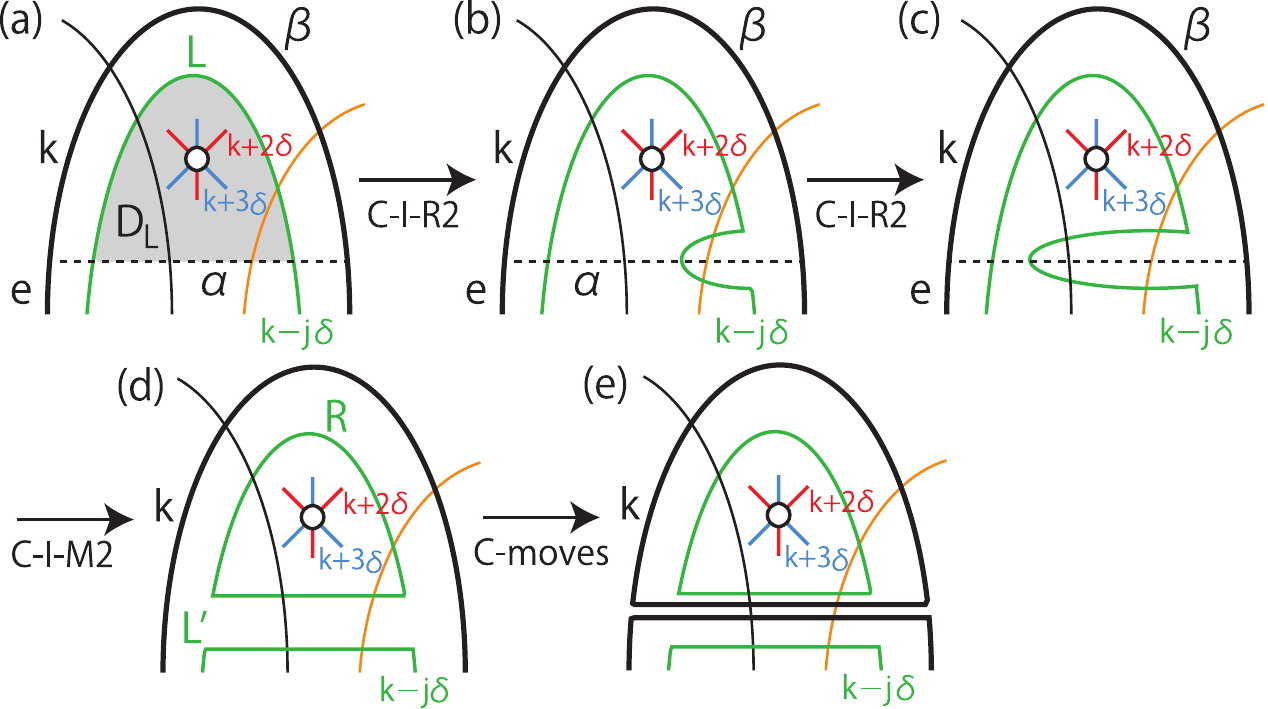}}
\caption{\label{fig13}
The gray region is the disk $D_L$,
$k$ is a label, $\delta\in\{+1,-1\}$,
$j$ is a positive integer.}
\end{figure}

\begin{lemma}
\label{DiskLemmaWhiteVertex}
{$($Disk Lemma with white vertices$)$}
Let $\Gamma$ be a chart, $k$ a label of $\Gamma$.
Let $e$ be an internal edge or a ring or a hoop of label $k$.
Let $D$ be a disk with a boundary arc pair $(\alpha,\beta)$
with $\Gamma_k\cap\partial D=\beta\subset e$ and
$\Gamma_{k+\delta}\cap\partial D=\emptyset$
for some $\delta\in\{+1,-1\}$.
Suppose that if an edge of $\Gamma$ intersects ${\rm Int}\alpha$,
then the edge transversely intersects the arc $\alpha$
$($see Fig.~\ref{fig12}$($a$))$.
Let $V$ be a neighborhood of $\alpha$.
If $D$ does not contain any white vertices in
 $\Gamma_{k+\delta}\cup(\cup_{i=0}^{\infty}\Gamma_{k-i\delta})$,
then we can replace the edge $e$ by the set $(e-\beta)\cup\alpha$
by C-moves in $V$ keeping 
$\cup_{i=1}^{\infty}\Gamma_{k+i\delta}$ fixed without increasing
the complexity of $\Gamma$
$($see Fig.~\ref{fig12}$($b$))$.
\end{lemma}

\begin{Proof}
Since $D$ does not contain any white vertex $\Gamma_{k+\delta}\cup(\cup_{i=0}^{\infty}\Gamma_{k-i\delta})$,
the disk $D$ does not contain any black vertices in $\Gamma_{k+\delta}\cup(\cup_{i=0}^{\infty}\Gamma_{k-i\delta})$.
Moreover, 
\begin{enumerate}
\item[(1)] $\Gamma_{k-i\delta}\cap D$ consists of proper arcs
of $D$ for all $i\geqq0$. 
\end{enumerate}

First, we shall show that
by applying C-moves in $V$,
we can assume that there is no $(D,\alpha)$-arcs of label $k-i\delta$ for all $i>0$.
We prove by induction on the number of $(D,\alpha)$-arcs of 
label $k-i\delta$ for all $i>0$.
Let $n$ be the number of  $(D,\alpha)$-arcs of 
label $k-i\delta$ for all $i>0$.

Suppose $n>0$.
Then there exists a $(D,\alpha)$-arc $L$ of label $k-j\delta$
for some $j>0$ (see Fig.~\ref{fig13}(a))
such that the disk $D_L$ with a boundary arc pair $(L,L_\alpha)$
contains no other $(D,\alpha)$-arcs of label $k-i\delta$ 
for $i>0$,
where $L_\alpha$ is an arc in $\alpha$.
In particular,
the condition $\Gamma_k\cap \partial D=\beta$ 
(i.e. $\Gamma_k\cap{\rm Int}\alpha=\emptyset$)
implies that 
$(\Gamma_{k-(j-1)\delta}\cup \Gamma_{k-j\delta}\cup\Gamma_{k-(j+1)\delta})\cap{\rm Int}L_\alpha=\emptyset$.

Let $\widetilde{L}$ be the connected component of 
$\Gamma_{k-j\delta}\cap(D\cup V)$ containing the arc $L$.
Then by deforming $\widetilde{L}$ in $V$ by
C-I-R2 moves,
we can push an end point of $L$ near the other end point of $L$
along $L_\alpha$ (see Fig.~\ref{fig13}(b),(c))
so that 
we can assume $\Gamma\cap{\rm Int}L_\alpha=\emptyset$.
By applying a C-I-M2 move (see Fig.~\ref{fig13}(d)),
we can split the arc  $\widetilde{L}$ 
to a ring (or a hoop) $R$ and an arc $L'$ 
to get a new chart $\Gamma'$ with 
$(R\cup L')\cap\alpha=\emptyset$.
Hence by induction,
we can assume that the chart does not contain 
$(D,\alpha)$-arcs of label $k-i\delta$ for all $i>0$.

Thus by (1), we can assume that
$\Gamma_{k-\delta}\cap \alpha=\emptyset$.
Hence the two conditions
$\Gamma_k\cap\partial D=\beta$ and
$\Gamma_{k+\delta}\cap \partial D=\emptyset$
imply 
$(\Gamma_{k-\delta}\cup\Gamma_k\cup\Gamma_{k+\delta})\cap{\rm Int}\alpha=\emptyset$.
Similarly, we can deform $\beta$ by C-I-R2 moves and a C-I-M2 move in $V$, and
we can replace the edge $e$ by $(e-\beta)\cup\alpha$
(see Fig.~\ref{fig13}(e)).
\end{Proof}

\begin{figure}[htb]
\centerline{\includegraphics{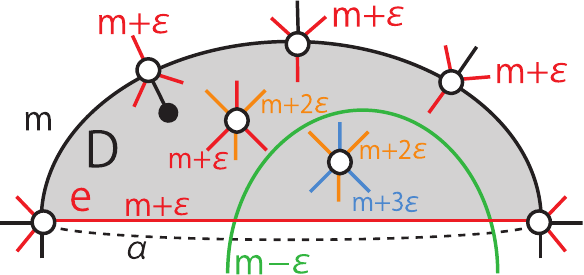}}
\caption{\label{fig14}
The gray region is the disk $D$,
$m$ is a label, $\varepsilon\in\{+1,-1\}$. }
\end{figure}

\begin{corollary}
\label{CorDiskLemmaWhiteVertice}
$($Corollary of Disk Lemma with white vertices$)$
Let $\Gamma$ be a chart, $m$ a label of $\Gamma$.
Let $D$ be a disk with a boundary arc pair $(e,\beta)$
such that $e$ is an internal edge of label $m+\varepsilon$
for some $\varepsilon\in\{+1,-1\}$,
$\beta\subset\Gamma_m$ and
$\Gamma_{m-\varepsilon}\cap\beta=\emptyset$
$($see Fig.~\ref{fig14}$)$.
Suppose that
 ${\rm Int}D$ does not contain any white vertices
in $\cup^{\infty}_{i=0}\Gamma_{m-i\varepsilon}$.
Then for a neighborhood $V$ of $e$,
there exists a chart $\Gamma'$ obtained by C-moves in $V$
keeping $\cup^{\infty}_{i=0}\Gamma_{m+i\varepsilon}$ fixed
without increasing the complexity of $\Gamma$
such that $\Gamma'_{m-\varepsilon}\cap e=\emptyset$.
\end{corollary}

\begin{Proof}
Since $\Gamma_{m-\varepsilon}\cap\beta=\emptyset$,
the arc $\beta$ does not contain any white vertices in 
$\Gamma_{m-\varepsilon}$.
Moreover, since ${\rm Int}D$ does not contain any white vertices
 in $\Gamma_{m-\varepsilon}$,
\begin{enumerate}
\item[(1)] the disk $D$ does not contain any white vertices 
in $\Gamma_{m-\varepsilon}$. 
\end{enumerate}

Let $\alpha$ be a simple arc parallel to $e$
with $\partial \alpha=\partial e$ and
such that $\alpha\cup\beta$ bounds a disk $D'$ containing
the disk $D$
(see Fig.~\ref{fig14}).
We can assume that if an edge of $\Gamma$ intersects 
${\rm Int}\alpha$,
then the edge transversely intersects the arc $\alpha$.

We prove this corollary by induction on the number of
points in $\Gamma_{m-\varepsilon}\cap e$.
Suppose that $\Gamma_{m-\varepsilon}\cap e\not=\emptyset$
(i.e. $\Gamma_{m-\varepsilon}\cap\alpha\not=\emptyset$).
Then by (1),
there exists a $(D',\alpha)$-arc $L$ of label $m-\varepsilon$
such that the disk $D_L$ with a boundary arc pair
$(L,L_\alpha)$ does not contain any other 
$(D',\alpha)$-arc of label $m-\varepsilon$,
where $L_\alpha$ is an arc in $\alpha$.
Hence 
\begin{enumerate}
\item[(2)] $\Gamma_{m-\varepsilon}\cap\partial D_L=L$
and
$\Gamma_m\cap\partial D_L=\emptyset$.
\end{enumerate}
Since ${\rm Int}D$ does not contain any white vertices
in $\cup^{\infty}_{i=0}\Gamma_{m-i\varepsilon}$ 
by the condition of this lemma,
the disk $D_L$ does not contain
any white vertices in 
$\Gamma_m\cup(\cup^{\infty}_{i=0}\Gamma_{(m-\varepsilon)-i\varepsilon})$.
Thus by (2) and 
Lemma~\ref{DiskLemmaWhiteVertex}(Disk Lemma with white vertices),
we obtain a chart $\Gamma'$ by
moving the arc $L$ of label $m-\varepsilon$ to $L_\alpha$
by C-moves keeping $\cup^{\infty}_{i=0}\Gamma_{m+i\varepsilon}$
fixed 
so that the number of points in $\Gamma_{m-\varepsilon}'\cap e$
is less than the number of 
points in $\Gamma_{m-\varepsilon}\cap e$.
Hence by induction,
we obtain a desired chart $\Gamma''$
with $\Gamma_{m-\varepsilon}''\cap e=\emptyset$.
\end{Proof}

\begin{figure}[ht]
\centerline{\includegraphics{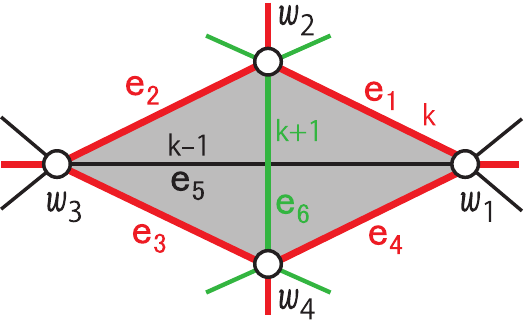}}
\vspace{5mm}
\caption{\label{fig15}
The gray region is the M4-disk. }
\end{figure}

Let $\Gamma$ be a chart and 
$k$ a label of $\Gamma$.
If a disk $D$ satisfies the following three conditions,
then $D$ is called an 
{\it M4-disk of label $k$}
 (see Fig.~\ref{fig15}).
\begin{enumerate}
\item[(i)] 
$\partial D$ consists of four internal edges 
$e_1,e_2,e_3,e_4$ of label $k$ 
situated on $\partial D$ 
in this order.
\item[(ii)] 
Set $w_1=e_1\cap e_4,w_2=e_1\cap e_2,
w_3=e_2\cap e_3,w_4=e_3\cap e_4$.
Then
\begin{enumerate}
\item[(a)] 
$D\cap\Gamma_{k-1}$ consists of an internal edge $e_5$ 
connecting $w_1$ and $w_3$, 
and
\item[(b)] 
$D\cap\Gamma_{k+1}$ consists of an internal edge $e_6$ 
connecting $w_2$ and $w_4$. 
\end{enumerate}
\item[(iii)] 
Int$D$ does not contain any white vertex.
\end{enumerate}
We call the union 
$X=\cup_{i=1}^6 e_i$ 
the {\it M4-pseudo chart} for the disk $D$.

\begin{lemma}
\label{M4-diskLemma}
{\rm (\cite[Lemma 7.3]{Chart33})}
Let $\Gamma$ be a chart, and 
$k$ a label of $\Gamma$. 
Suppose that $D$ is an M4-disk 
of label $k$ 
with an M4-pseudo chart $X$. 
Then by deforming $\Gamma$ 
in a regular neighbourhood of $D$ 
without increasing the complexity of $\Gamma$, 
the chart $\Gamma$ is 
C-move equivalent to a chart $\Gamma'$ 
with
$D\cap(\cup_{i=k-2}^{k+2}\Gamma'_i)=X$.
\end{lemma}

\begin{figure}[thb]
\centerline{\includegraphics{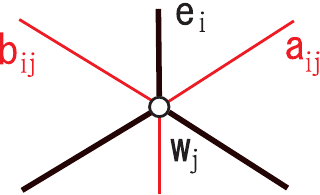}}
\caption{\label{fig16}
The three edges $a_{ij},e_i,b_{ij}$ are consecutive edges around the white vertex $w_j$.}
\end{figure}

In our argument,
we often need a name for an unnamed edge by using a given edge and a given white vertex.
For the convenience,
we use the following naming:
Let $e',e_i,e''$ be three consecutive edges containing  a white vertex $w_j$. Here, 
the two edges $e'$ and $e''$ are unnamed edges. 
There are six arcs in a neighborhood $U$ of the white vertex $w_j$. 
If the three arcs $e'\cap U$, $e_i \cap U$, $e'' \cap U$ lie anticlockwise around the white vertex $w_j$ in this order, 
then $e'$ and $e''$ are denoted by $a_{ij}$ and $b_{ij}$ 
respectively (see Fig.~\ref{fig16}).
There is a possibility $a_{ij}=b_{ij}$ if they are contained in a loop.

\begin{figure}[htb]
\centerline{\includegraphics{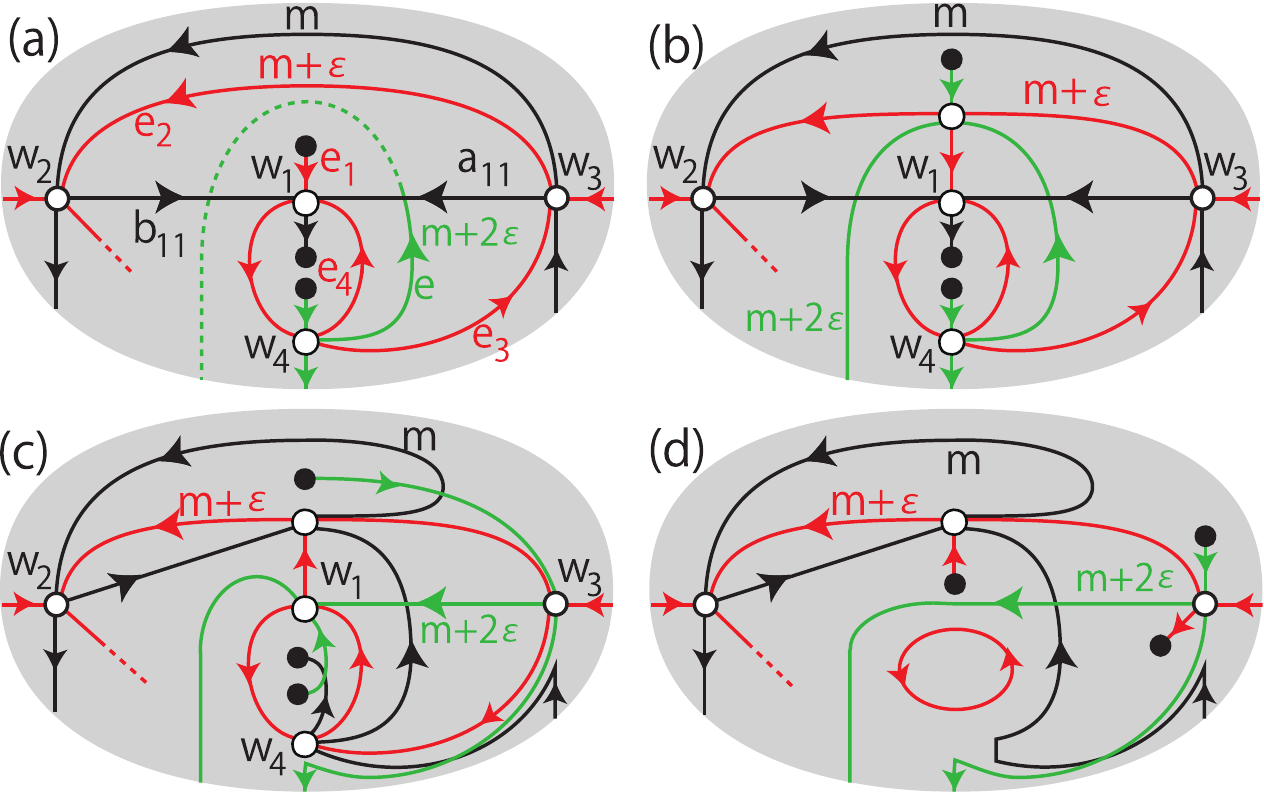}}
\caption{\label{fig17}
The gray regions are the disk $D$,
$m$ is a label, $\varepsilon\in\{+1,-1\}$.}
\end{figure}

\begin{lemma}
\label{ApplicationCIR4Move}
Let $\Gamma$ be a chart, and $m$ a label of $\Gamma$.
If $\Gamma$ contains the pseudo chart in a disk $D$
as shown in Fig.~\ref{fig17}$($a$)$,
and if $w(\Gamma\cap D)=4$,
then $\Gamma$ is not a minimal chart. 
\end{lemma}

\begin{Proof}
Suppose that $\Gamma$ is minimal. 
We use the notations as shown in 
Fig.~\ref{fig17}(a),
where $e_1$ is a terminal edge of label $m+\varepsilon$ at $w_1$,
$e_2,e_3,e_4$ are internal edges of label $m+\varepsilon$
with $\partial e_2=\{w_2,w_3\}$,
$\partial e_3=\{w_3,w_4\}$,
$\partial e_4=\{w_1,w_4\}$,
and $e$ is an internal edge of label $m+2\varepsilon$
at $w_4$.
Let $a_{11},b_{11}$ be internal edges of label $m$
with $\partial a_{11}=\{w_1,w_3\}$,
$\partial b_{11}=\{w_1,w_2\}$.
Let $E_1,E_2$ be disks in $D$
with $\partial E_1=a_{11}\cup e_3\cup e_4$
and $\partial E_2=a_{11}\cup b_{11}\cup e_2$.
Since $w(\Gamma\cap D)=4$ by the condition of this lemma,
\begin{enumerate}
\item[(1)] neither ${\rm Int}E_1$ nor ${\rm Int}E_2$ 
contains white vertices.
\end{enumerate}

Let $\alpha$ be an arc connecting the black vertex in $e_1$
and a point in $e_2$ with $\alpha\subset E_2$.

{\bf Claim. } $(\Gamma_m\cup\Gamma_{m+\varepsilon}\cup\Gamma_{m+2\varepsilon})\cap{\rm Int}\alpha=e\cap{\rm Int}\alpha=$one point
by C-moves in $D$ without increasing the complexity of $\Gamma$.

{\it Proof of Claim.} 
Apply New Disk Lemma(Lemma~\ref{NewDiskLemma}) 
for the disk $E_1$,
we can assume that the chart $\Gamma$ is $(E_1,a_{11})$-arc free.
Thus $\Gamma_{m+2\varepsilon}\cap a_{11}=e\cap a_{11}=$~one point.
Because, if $\Gamma_{m+2\varepsilon}\cap a_{11}$ consists at least two points,
then by (1)
there exists a proper arc $\gamma$ of $E_1$ in an internal edge
or a ring of label $m+2\varepsilon$.
Since 
$e_3\cup e_4\subset \Gamma_{m+\varepsilon}\cap \partial E_1$,
we have $\partial \gamma\subset a_{11}$.
Hence the arc $\gamma$ is a $(E_1,a_{11})$-arc of label $m+2\varepsilon$.
This contradicts the fact 
that the chart $\Gamma$ is $(E_1,a_{11})$-arc free. Thus
\begin{enumerate}
\item[(2)] $\Gamma_{m+2\varepsilon}\cap a_{11}=e\cap a_{11}=$ one point.
\end{enumerate}

Let $N(e_1)$ be a regular neighborhood of the terminal edge $e_1$
in $E_2$. 
Set $E_2'=Cl(E_1-N(e_1))$ and $b_{11}'=b_{11}\cap E_2'$.
Then by (1),
the disk $E_2'$ does not contain any black vertices.
Apply New Disk Lemma(Lemma~\ref{NewDiskLemma}) 
for the disk $E_2'$,
we can assume that 
the chart $\Gamma$ is $(E_2',b_{11}')$-arc free.
Hence the chart $\Gamma$ is $(E_2,b_{11})$-arc free.
Thus by the similar way as above,
we can show that 
$\Gamma_{m+2\varepsilon}\cap b_{11}=e\cap b_{11}=$one point.
Hence
\begin{enumerate}
\item[(3)] $\Gamma_{m+2\varepsilon}\cap E_2$ is one proper arc of $E_2$.
\end{enumerate}

Since $\partial E_2\subset \Gamma_m\cup\Gamma_{m+\varepsilon}$
and since $\Gamma$ is minimal,
the disk $E_2$ does not intersect any ring of label $m$ or $m+\varepsilon$
by (1) and Assumption~\ref{Ring}.
Hence $(\Gamma_m\cup\Gamma_{m+\varepsilon})\cap E_2=e_1\cup e_2\cup a_{11}\cup b_{11}$.
Thus $(\Gamma_m\cup\Gamma_{m+\varepsilon})\cap{\rm Int}\alpha=\emptyset$.
Hence by (2) and (3), we have
$(\Gamma_m\cup\Gamma_{m+\varepsilon}\cup\Gamma_{m+2\varepsilon})\cap{\rm Int}\alpha=e\cap{\rm Int}\alpha=$one point.
Thus Claim holds. {\hfill {$\square$}\vspace{1.5em}}

Hence by C-II moves and C-I-R2 moves,
we can assume that $\Gamma\cap {\rm Int}\alpha=e\cap {\rm Int}\alpha=$one point.
Thus we can apply a C-III move 
among the three edges $e_1,e,e_2$, and 
we obtain the pseudo chart as shown in Fig.~\ref{fig17}(b).
Then we can apply a C-I-R4 move by Lemma~\ref{M4-diskLemma},
and we obtain the pseudo chart as shown in
 Fig.~\ref{fig17}(c).
Hence we obtain a terminal edge of label $m$ at $w_4$ and
a terminal edge of label $m+2\varepsilon$ at $w_1$ such that 
neither two terminal edges are middle at $w_1$ or $w_4$.
Thus by C-III moves,
the number of white vertices decreases
(see Fig.~\ref{fig17}(d)).
This is a contradiction.
Therefore $\Gamma$ is not minimal.
We complete the proof of Lemma~\ref{ApplicationCIR4Move}.
\end{Proof}



\section{IO-Calculation}

\label{s:IOC}

In this section,
we review IO-Calculation.

Let $\Gamma$ be a chart,
 and $v$ a vertex. 
Let $\alpha$ be a short arc of $\Gamma$ in a small neighborhood of $v$ such that $v$ is an endpoint of $\alpha$. 
If the arc $\alpha$ is oriented to $v$, then $\alpha$ is called {\it an inward arc}, 
and otherwise $\alpha$ is called {\it an outward arc}.

\begin{figure}
\centerline{\includegraphics{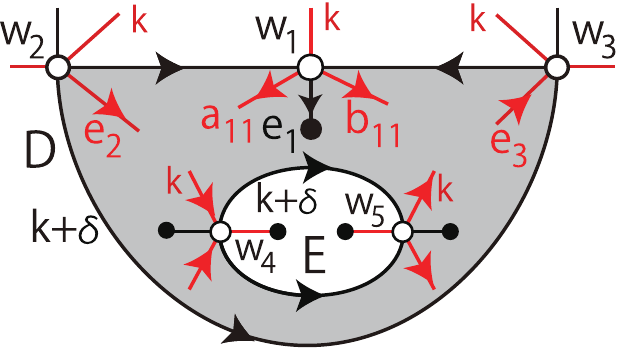}}
\caption{\label{fig18} The gray region is the region $F$, $k$ is a label, $\delta\in\{+1,-1\}$.}
\end{figure}

Let $\Gamma$ be an $n$-chart. 
Let $F$ be a closed domain with $\partial F\subset \Gamma_{k-1}\cup\Gamma_{k}\cup \Gamma_{k+1}$ for some label $k$ of $\Gamma$, where $\Gamma_0=\emptyset$ and $\Gamma_{n}=\emptyset$. 
By Condition (iii) for charts,
in a small neighborhood of each white vertex, there are three inward arcs and three outward arcs.
Also in a small neighborhood of each black vertex, there exists only one inward arc or one outward arc.
We often use the following fact, 
when we fix (inward or outward) arcs 
near white vertices and black vertices: 
\begin{enumerate}
\item[$(*)$]
{\it The number of inward arcs contained in $F\cap \Gamma_k$ is equal to the number of outward arcs in $F\cap \Gamma_k$.
}
\end{enumerate}
When we use this fact, 
we say that we use {\it IO-Calculation with respect to $\Gamma_k$ in $F$}.
For example, in a minimal chart $\Gamma$, 
consider the pseudo chart as shown in Fig.~\ref{fig18} 
where
\begin{enumerate}
\item[(1)] $D$ is a $3$-angled disk of $\Gamma_{k+\delta}$ 
with one feeler $e_1$ for some $\delta\in\{+1,-1\}$,
\item[(2)] $E$ is a $2$-angled disk of $\Gamma_{k+\delta}$ 
without feelers in $D$ with $F=Cl(D-E)$,
\item[(3)] $a_{11},b_{11},e_2$ are internal edges (possibly terminal edges) of label $k$ oriented outward at $w_1,w_1,w_2$, respectively, \item[(4)]  $e_3$ is an internal edge (possibly a terminal edge) 
of label $k$ oriented inward at $w_3$, 
\item[(5)] none of $a_{11},b_{11},e_3$ are
 middle at $w_1$ or $w_3$.
\end{enumerate}
Then we can show that $w(\Gamma\cap{\rm Int}F)\ge1$.
Suppose $w(\Gamma\cap{\rm Int}F)=0$.
Let $e_4,e_4'$ be
internal edges (possibly terminal edges)
of label $k$ oriented inward at $w_4$,
and $e_5,e_5'$
internal edges (possibly terminal edges)
of label $k$ oriented outward at $w_5$,
Then by (5) and Assumption~\ref{AssumeTerminal},
\begin{enumerate}
\item[(6)] none of $a_{11},b_{11},e_3,e_4,e_4',e_5,e_5'$ are 
terminal edges. 
\end{enumerate}

If the edge $e_2$ is a terminal edge,
then 
by (3),(4) and (6) 
the number of inward arcs in $F\cap \Gamma_k$ is four,  
but the number of outward arcs in $F\cap \Gamma_k$ is five. 
This contradicts the fact $(*)$. 
If $e_2$ is not a terminal edge,
then
by (3),(4) and (6) 
the number of inward arcs in $F\cap \Gamma_k$ is three,  
but the number of outward arcs in $F\cap \Gamma_k$ is five. 
This contradicts the fact~$(*)$. 
Thus $w(\Gamma\cap{\rm Int}F)\ge1$.
Instead of the above argument, 
\begin{enumerate}
\item[]
{\it we have $w(\Gamma\cap{\rm Int}F)\ge1$ 
by IO-Calculation with respect to $\Gamma_{k}$ in $F$.}
\end{enumerate}

\section{Lenses}
\label{s:Lens}

In this section, we review a useful lemma for a disk called a lens.

\begin{figure}[htb]
\centerline{\includegraphics{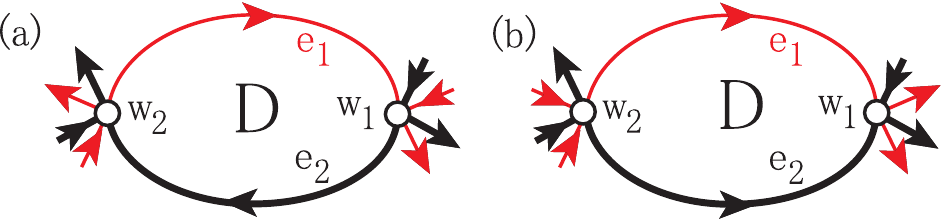}}
\caption{\label{fig19}
Lenses.}
\end{figure}

Let $\Gamma$ be a chart. 
Let $D$ be a disk 
such that 
\begin{enumerate}
\item[(1)] the boundary $\partial D$ consists of an internal edge $e_1$ of label $m$ and an internal edge $e_2$ of label ${m+1}$, and 
\item[(2)] any edge containing a white vertex in $e_1$ does not intersect the open disk Int$D$.
\end{enumerate}
Note that $\partial D$ may contain crossings.
Let $w_1$ and $w_2$ be the white vertices in $e_1$. 
If the disk $D$ satisfies one of the following conditions, then $D$ is called  {\it a lens of type $(m,m+1)$}
(see Fig.~\ref{fig19}):
\begin{enumerate}
	\item[(i)] Neither $e_1$ nor $e_2$ contains a middle arc. 
	\item[(ii)] One of the two edges $e_1$ and $e_2$ contains middle arcs at both white vertices $w_1$ and $w_2$ simultaneously.
\end{enumerate}

\begin{lemma}{\rm (\cite[Corollary 1.3]{ChartAppII})}
\label{NoLens}
 There is no lens in any minimal chart with 
at most seven white vertices.
\end{lemma}

\begin{figure}
\centerline{\includegraphics{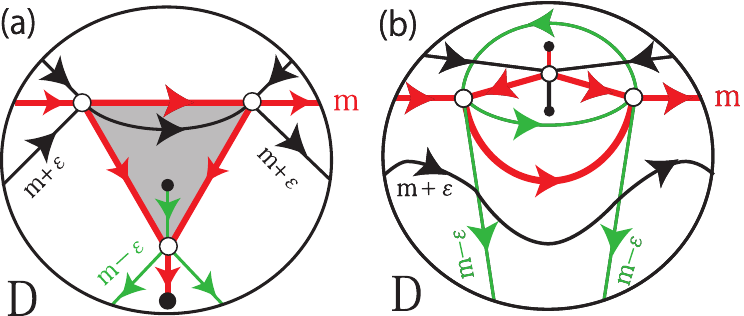}}
\caption{\label{fig20}
The gray region is the 3-angled disk $D_1$. 
The thick lines are edges of label $m$,
and $\varepsilon\in\{+1,-1\}$.}
\end{figure}

\begin{lemma}
\label{CorTriangle}
{\rm (\cite[Corollary 13.4]{ChartAppX})}
For a chart $\Gamma$,
if there exists  a $3$-angled disk $D_1$ of $\Gamma_m$ without feelers in a disk $D$ as shown in Fig.~\ref{fig20}$($a$)$ and if $w(\Gamma\cap${\rm Int}$D_1)=0$,
then there exists a chart obtained from $\Gamma$ by C-moves in $D$ which contains the pseudo chart in $D$ as shown in Fig.~\ref{fig20}$($b$)$. 
\end{lemma}

\begin{lemma}
{\rm (\cite[Theorem 1.1]{ChartAppIX})}
\label{NoMinimal(43)} 
There is 
no minimal chart of 
type $(4,3)$.
\end{lemma}

\begin{lemma}
\label{LemmaWithTerminal3}
{\rm (\cite[Lemma 3.2(1)]{ChartAppV})}
Let $\Gamma$ be a minimal chart,
and $m$ a label of $\Gamma$.
Let $G$ be a connected component of $\Gamma_m$.
If $1\le w(G)$, then $2\le w(G)$.
\end{lemma}



\section{Case of the oval}
\label{s:Oval}

In this section, we show that 
for any minimal chart $\Gamma$ of type $(m;5,2)$,
the graph $\Gamma_m$ does not contain an oval
as shown in Fig.~\ref{fig09}(b).

Let $\Gamma$ be a minimal chart of type $(m;5,2)$.
Suppose that $\Gamma_m$ contains an oval $G_1$.
Then by 
Lemma~\ref{GammaMFiveWhite}(c) and Lemma~\ref{LemmaNoLoop},
the graph $\Gamma_m$ contains a skew $\theta$-curve $G_2$.
The graph $G_2$ divides $S^2$ into three disks.
One of the three disks is a 2-angled disk, say $D_1$.
One of the three disks is a 3-angled disk with one feeler $e_1$,
say $D_2$.
Let $D_3$ be the third disk.
Since $D_2$ has exactly one feeler $e_1$,
by Lemma~\ref{Theorem3AngledDisk} we have
\begin{enumerate}
\item[(a)] $w(\Gamma\cap{\rm Int}D_2)\geqq1$.
\end{enumerate}

\begin{figure}[htb]
\centerline{\includegraphics{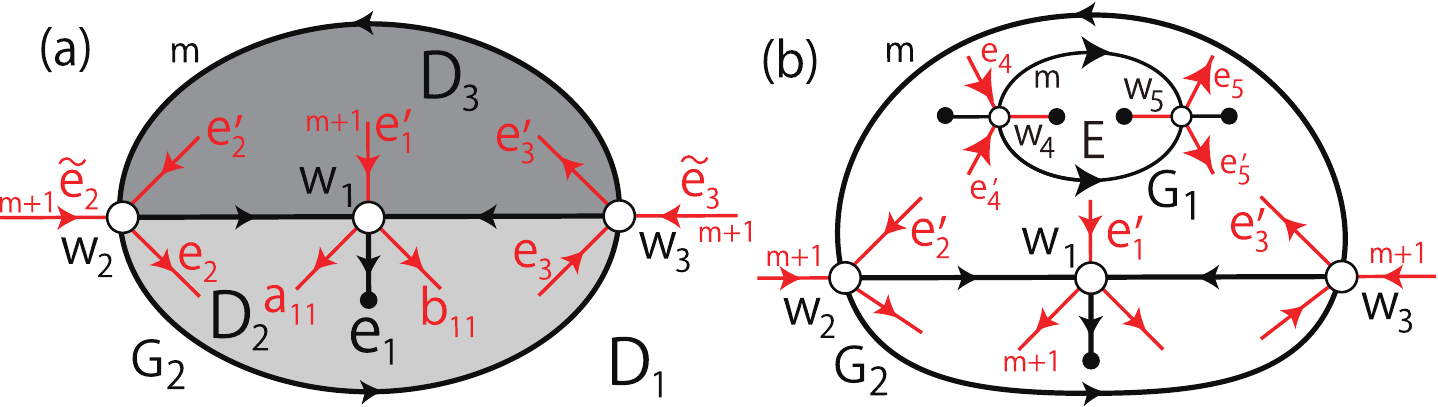}}
\caption{\label{fig21}
The light gray region is the disk $D_2$.
The dark gray region is the disk $D_3$.
}
\end{figure}

Without loss of generality,
we can assume that the terminal edge $e_1$
is oriented outward at the white vertex $w_1$ in $e_1$.
Since $e_1$ is middle at $w_1$ 
by Assumption~\ref{AssumeTerminal},
the two internal edges of label $m$ at $w_1$
are oriented inward at $w_1$.
Hence by Lemma~\ref{Theorem2AngledDisk}
\begin{enumerate}
\item[(b)] $w(\Gamma\cap{\rm Int}D_1)\geqq1$.
\end{enumerate}
Let $w_2,w_3$ be the white vertices in 
the skew $\theta$-curve $G_2$ different from $w_1$.
Without loss of generality,
we can assume that 
the intersection $D_1\cap D_2$ is oriented from $w_2$ to $w_3$.
By looking around $w_2$,
the edge $D_1\cap D_3$ is oriented from $w_3$ to $w_2$.
Therefore, 
the chart $\Gamma$ contains the pseudo chart as shown
in Fig.~\ref{fig21}(a).
From now on throughout this section,
we use the notations as shown in Fig.~\ref{fig21}(a).

\begin{lemma}
\label{NotContainD3}
Let $\Gamma$ be a minimal chart of type $(m;5,2)$.
If $\Gamma_m$ contains an oval $G_1$ and 
a skew $\theta$-curve $G_2$,
then $G_1$ is not contained in the $3$-angled disk $D_3$
without feelers.
\end{lemma}

\begin{Proof}
Suppose $G_1\subset D_3$.
Then by Conditions (a) and (b) of this section,
the condition $w(\Gamma)=7$ implies that
$w(\Gamma\cap{\rm Int}D_3)=2$.

Let $E$ be the 2-angled disk of $\Gamma_m$
in $D_3$ with $\partial E\subset G_1$.
Then the condition $w(\Gamma\cap{\rm Int}D_3)=2$
implies that 
$w(\Gamma\cap{\rm Int}E)=0$.
Thus by Lemma~\ref{Theorem2AngledDisk},
a regular neighborhood of $E$ contains the pseudo chart
as shown in Fig.~\ref{fig03}(b).
Hence $\Gamma$ contains the pseudo chart as shown
in Fig.~\ref{fig21}(b),
where
\begin{enumerate}
\item[(1)] $e_1',e_2',e_4,e_4'$ are internal edges 
(possibly terminal edges) of label $m+1$
oriented inward at
$w_1,w_2,w_4,w_4$, respectively.
\end{enumerate}
Moreover, none of $e_2',e_4,e_4'$ are middle at $w_2$ or $w_4$.
Thus by Assumption~\ref{AssumeTerminal},
\begin{enumerate}
\item[(2)] none of $e_2',e_4,e_4'$ are terminal edges. 
\end{enumerate}
Hence the condition $w(\Gamma\cap{\rm Int}D_3)=2$ implies that
\begin{enumerate}
\item[(3)] the edge $e_1'$ must be a terminal edge,
\item[(4)] none of $e_3',e_5,e_5'$ are terminal edges.
\end{enumerate}

For the edge $e_2'$,
there are three cases:
(i) $e_2'=e_3'$,
(ii) $e_2'=e_5$,
(iii) $e_2'=e_5'$.

{\bf Case (i).}
Since $e_2'=e_3'$,
we have $e_4=e_5$ and $e_4'=e_5'$.
Thus there exist two lenses in $D_3$.
This contradicts Lemma~\ref{NoLens}.
Hence Case (i) does not occur.

{\bf Case (ii).}
Since $e_2'=e_5$,
we have $e_4=e_3'$ and $e_4'=e_5'$.
Thus there exists a lens in $D_3$.
This contradicts Lemma~\ref{NoLens}.
Hence Case (ii) does not occur.

{\bf Case (iii).}
Since $e_2'=e_5'$,
we have $e_4=e_5$ and $e_4'=e_3'$.
Thus there exists a lens in $D_3$.
This contradicts Lemma~\ref{NoLens}.
Hence Case (iii) does not occur.

Therefore, all the three cases do not occur.
Hence $G_1\not\subset D_3$.
\end{Proof}



\begin{lemma}
\label{NotContainD2}
Let $\Gamma$ be a minimal chart of type $(m;5,2)$.
If $\Gamma_m$ contains an oval $G_1$ and 
a skew $\theta$-curve $G_2$,
then $G_1$ is not contained in the $3$-angled disk $D_2$
with one feeler.
\end{lemma}

\begin{figure}[htb]
\centerline{\includegraphics{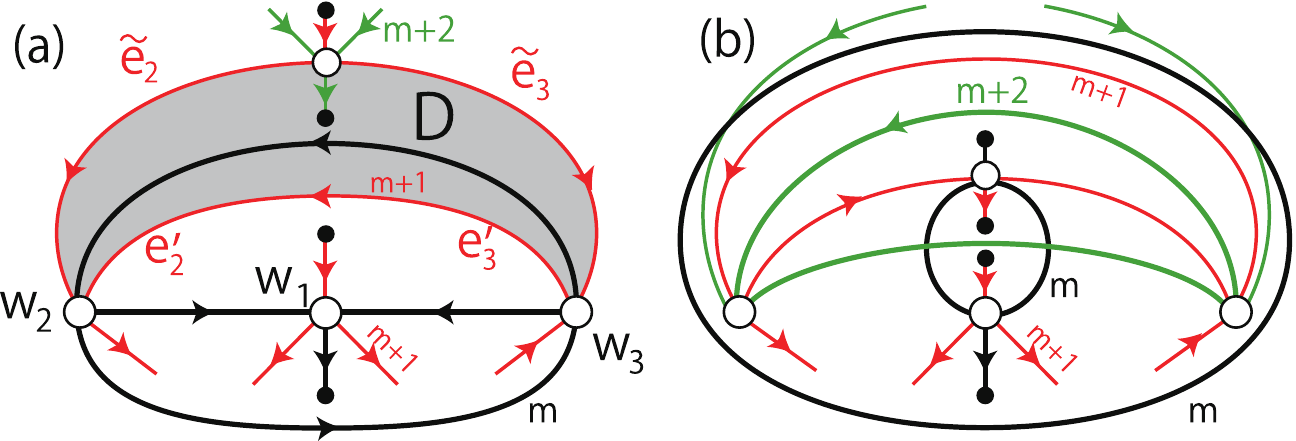}}
\caption{\label{fig22}
The gray region is the disk $D$.}
\end{figure}

\begin{Proof}
Suppose $G_1\subset D_2$.
We use the notations as shown in Fig.~\ref{fig21}(a).

{\bf Claim.} $w(\Gamma\cap {\rm Int}D_2)\geqq3$.

{\it Proof of Claim.}
 Let $E$ be the 2-angled disk of $\Gamma_m$ in $D_2$
with $\partial E\subset G_1$.
If $w(\Gamma\cap{\rm Int}E)\geqq1$,
then we have $w(\Gamma\cap{\rm Int}D_2)\geqq3$.

Now, suppose that $w(\Gamma\cap{\rm Int}E)=0$.
Then by Lemma~\ref{Theorem2AngledDisk},
a regular neighborhood of $E$ contains
the pseudo chart as shown in Fig.~\ref{fig03}(b).
Thus the 2-angled disk $E$ has no feelers.
Hence the chart $\Gamma$ contains
the pseudo chart as shown in Fig.~\ref{fig18}
where $k=m+1$ and $\delta=-1$.
Thus,
we have $w(\Gamma\cap({\rm Int}D_2-E))\geqq1$
by considering as $F=Cl(D_2-E)$  
in the example of IO-Calculation in 
Section~\ref{s:IOC}.
Hence we have $w(\Gamma\cap{\rm Int}D_2)\geqq3$.
Thus Claim holds.{\hfill {$\square$}\vspace{1.5em}}

By Claim and Condition~(b) of this section,
the condition $w(\Gamma)=7$ implies that
\begin{enumerate}
\item[(1)]
$w(\Gamma\cap {\rm Int}D_1)=1$,
 $w(\Gamma\cap {\rm Int}D_2)=3$, 
and 
$w(\Gamma\cap {\rm Int}D_3)=0$.
\end{enumerate}
Thus by Lemma~\ref{Theorem2AngledDisk2},
a regular neighborhood of $D_1$
contains one of RO-families of the two pseudo charts
as shown in Fig.~\ref{fig07}(b),(c).
Moreover,
by Lemma~\ref{Theorem3AngledDisk},
a regular neighborhood of $D_3$
contains one of the RO-family of the pseudo chart
as shown in Fig.~\ref{fig04}(b).
Hence $e_2'=e_3'$ and $\widetilde{e}_2\cap \widetilde{e}_3$
is a white vertex in $\Gamma_{m+1}\cap\Gamma_{m+2}$,
say $w_7$.
Let $e_7$ be the terminal edge of label $m+1$ at $w_7$, and
$D$ the 3-angled disk of $\Gamma_{m+1}$
in $D_1\cup D_3$
with $\partial D=e_2'\cup \widetilde{e}_2\cup \widetilde{e}_3$.
Then by (1), we have $w(\Gamma\cap {\rm Int}D)=0$.
Thus by Lemma~\ref{Theorem3AngledDisk},
a regular neighborhood of $D$
contains one of the RO-family of the pseudo chart
as shown in Fig.~\ref{fig04}(a).
Hence $e_7\not\subset D$ and
there exists a terminal edge of label $m+2$ at $w_7$
in $D$ (see Fig.~\ref{fig22}(a)).
We can apply Lemma~\ref{CorTriangle}
for the disk $D$.
Then we obtain the pseudo chart as shown in
Fig.~\ref{fig22}(b).
Thus we obtain a minimal chart of type $(m;4,3)$.
This contradicts Lemma~\ref{NoMinimal(43)}.
Hence $G_1\not\subset D_2$.
\end{Proof}

\begin{proposition}
\label{NoOvalSkewTheta}
Let $\Gamma$ be a minimal chart of type $(m;5,2)$.
Then $\Gamma_m$ does not contain an oval.
\end{proposition}

\begin{figure}[htb]
\centerline{\includegraphics{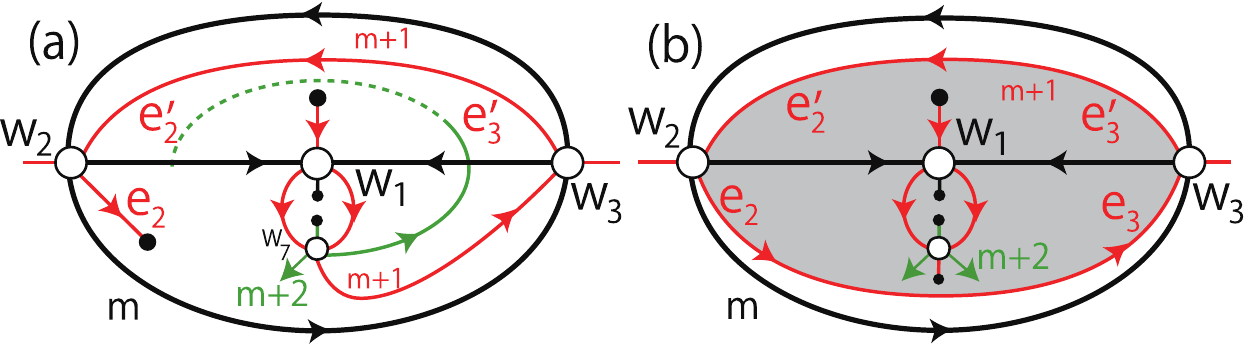}}
\caption{\label{fig23}
The gray region is the disk $D$.}
\end{figure}

\begin{Proof}
Suppose that $\Gamma_m$ contains an oval, say $G_1$.
Then by Lemma~\ref{GammaMFiveWhite}(c)
and Lemma~\ref{LemmaNoLoop},
the graph $\Gamma_m$ contains a skew $\theta$-curve, say $G_2$.
By the aurgument of the begining this section,
the graph $G_2$ is the graph as shown in 
Fig.~\ref{fig21}(a).
We use the notations as shown in Fig.~\ref{fig21}(a),
where 
\begin{enumerate}
\item[(1)] $\widetilde{e}_2,\widetilde{e}_3$ are internal edges
(possibly terminal edges) of label $m+1$ oriented inward at
$w_2,w_3$, respectively.
\end{enumerate}
Moreover, neither $\widetilde{e}_2$ nor $\widetilde{e}_3$ is
middle at $w_2$ or $w_3$.
Thus by Assumption~\ref{AssumeTerminal},
\begin{enumerate}
\item[(2)] neither $\widetilde{e}_2$ nor $\widetilde{e}_3$
is a terminal edge.
\end{enumerate}

By Lemma~\ref{NotContainD3} and Lemma~\ref{NotContainD2},
the oval $G_1$ is contained in the 2-angled disk $D_1$ of 
$\Gamma_m$ without feelers.

{\bf Claim.} $w(\Gamma\cap{\rm Int}D_1)\geqq3$.

{\it Proof of Claim.}
Let $E$ be the 2-angled disk of $\Gamma_m$ in $D_1$ with
$\partial E\subset G_1$.
If $w(\Gamma\cap{\rm Int}E)\geqq1$,
then $w(\Gamma\cap{\rm Int}D_1)\geqq3$.

Now suppose that 
$w(\Gamma\cap{\rm Int}E)=0$.
Thus by Lemma~\ref{Theorem2AngledDisk},
a regular neighborhood of $E$ contains
the pseudo chart 
as shown in Fig.~\ref{fig03}(b).
Let $w_4$ be the white vertex in $G_1$ and
$e_4$ the terminal edge in $G_1$
such that $e_4$ is oriented inward at the white vertex
$w_4$.
Let $a_{44},b_{44}$ be internal edges (possibly terminal edges)
of label $m+1$ oriented inward at $w_4$.
Then by Assumption~\ref{AssumeTerminal},
neither $a_{44}$ nor $b_{44}$
is middle at $w_4$.
Thus neither $a_{44}$ nor $b_{44}$
is a terminal edge.
Hence by (1) and (2),
we have $w(\Gamma\cap({\rm Int}D_1-E))\geqq1$
by IO-Calculation with respect to $\Gamma_{m+1}$
in $Cl(D_1-E)$.
Thus we have $w(\Gamma\cap{\rm Int}D_1)\geqq3$.
Therefore Claim holds.
{\hfill {$\square$}\vspace{1.5em}}

By Claim and Condition~(a) of this section,
the condition $w(\Gamma)=7$ implies that
\begin{enumerate}
\item[(3)] $w(\Gamma\cap{\rm Int}D_1)=3$,
 $w(\Gamma\cap{\rm Int}D_2)=1$ and 
$w(\Gamma\cap{\rm Int}D_3)=0$.
\end{enumerate}
Thus, by Lemma~\ref{Theorem3AngledDisk},
a regular neighborhood of $D_3$ contains
one of the RO-family of the pseudo chart
as shown in Fig.~\ref{fig04}(b).
Moreover, by Lemma~\ref{Theorem3AngledDiskChart9},
a regular neighborhood $N(D_2)$ of $D_2$ contains 
one of the RO-families of the two pseudo charts as
shown in Fig.~\ref{fig05}(e),(f).

Suppose that the neighborhood $N(D_2)$ contains
one of the RO-family of the pseudo chart as
shown in Fig.~\ref{fig05}(e)
(see Fig.~\ref{fig23}(a)).
Thus by (3),
the chart $\Gamma$ contains the pseudo chart as shown in 
Fig.~\ref{fig17}(a).
Hence by Lemma~\ref{ApplicationCIR4Move},
the chart $\Gamma$ is not minimal.
This is a contradiction.

Suppose that the neighborhood $N(D_2)$ contains
one of the RO-family of the pseudo chart as
shown in Fig.~\ref{fig05}(f) (see Fig.~\ref{fig23}(b)).
Then we have $e_2=e_3$.
Thus there exists a 2-angled disk $D$ of $\Gamma_{m+1}$
in $D_2\cup D_3$ with $\partial D=e_2\cup e_2'$.
Moreover, by (3),
the disk $D$ contains exactly two white vertices.
One of the two white vertices
is contained in $\Gamma_m\cap\Gamma_{m+1}$.
The other is contained in $\Gamma_{m+1}\cap\Gamma_{m+2}$.
Therefore there exists a connected component of $\Gamma_{m+2}$
with exactly one white vertex.
This contradicts Lemma~\ref{LemmaWithTerminal3}.

Hence we have a contradiction for the both cases. Therefore,
the graph $\Gamma_m$ does not contain an oval.
\end{Proof}

\section{Case of the graph as shown in Fig.~\ref{fig10}(h)}
\label{s:TypeH}

In this section,
we shall show that
for any minimal chart $\Gamma$ of type $(m;5,2)$,
the graph $\Gamma_m$ does not contain the graph as
shown in Fig.~\ref{fig10}(h).

\begin{figure}[htb]
\centerline{\includegraphics{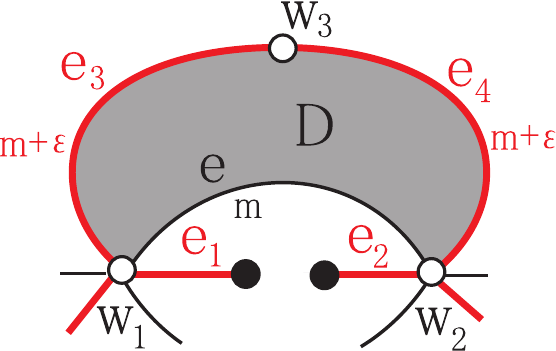}}
\caption{\label{fig24}
The white vertices $w_1$ and $w_2$ are in $\Gamma_m\cap\Gamma_{m+\varepsilon}$, and $\varepsilon\in\{+1,-1\}$.}
\end{figure}

\begin{lemma}
\label{CorDiskLemma}
{\rm (\cite[Lemma 5.4]{ChartApp1})}
If a minimal chart $\Gamma$ contains the pseudo chart 
as shown in Fig.~\ref{fig24}, 
then the interior of the disk $D$ contains at least one white vertex, 
where $D$ is the disk with the boundary $e_3\cup e_4\cup e$.
\end{lemma}

\begin{figure}
\begin{center}
\includegraphics{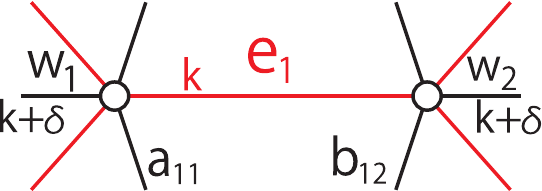}
\end{center}
\caption{\label{fig25} The edge $e_1$ is of label $k$,
and $\delta\in\{+1,-1\}.$ }
\end{figure}

\begin{lemma}
\label{C-I-M2A11B12}
{\rm (\cite[Lemma 3.3]{ChartAppX})}
Let $\Gamma$ be a chart, and
$k$ a label of $\Gamma$.
Let $e_1$ be an internal edge of label $k$
with two white vertices $w_1$ and $w_2$
$($see Fig.~\ref{fig25}$)$.
Suppose that $w_1,w_2\in\Gamma_{k+\delta}$
for some $\delta\in\{+1,-1\}$,
and suppose that one of the two edges $a_{11},b_{12}$
is a terminal edge.
If $\Gamma_{k+2\delta}\cap e_1=\emptyset$,
and if $\Gamma$ satisfies one of the following four conditions,
then $\Gamma$ is not a minimal chart. 
\begin{enumerate}
\item[{\rm (a)}]
The two edges $a_{11},b_{12}$ are oriented inward $($or outward$)$ at $w_1,w_2$, respectively.
\item[{\rm (b)}] 
The edge $a_{11}$ $($resp. $b_{12})$ is a terminal edge, and
$b_{12}$ $($resp. $a_{11})$ is not middle at the white vertex
different from $w_2$ $($resp. $w_{1})$.
\item[{\rm (c)}] 
The two edges $a_{11},b_{12}$ are middle at  $w_1,w_2$, respectively.
\item[{\rm (d)}] 
Both of  $a_{11},b_{12}$ are terminal edges.
\end{enumerate}
\end{lemma}

\begin{figure}[thb]
\centerline{\includegraphics{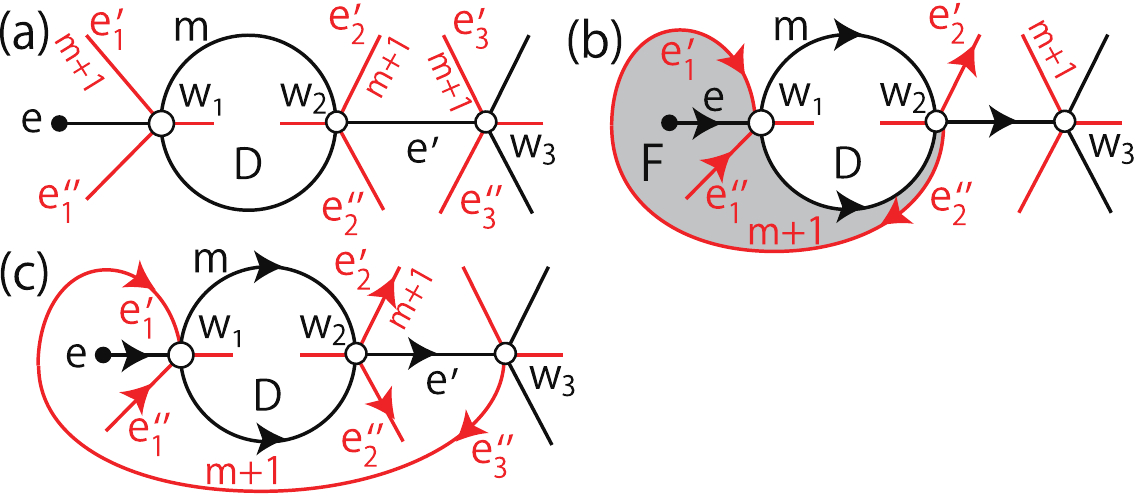}}
\caption{\label{fig26}
The gray region is the region $F$.}
\end{figure}

\begin{lemma}
\label{LemmaTypeGH}
Let $\Gamma$ be a minimal chart of type $(m;5,2)$.
Suppose that $\Gamma_m$ contains one of the two graphs
as shown in Fig.~\ref{fig10}$($g$)$,$($h$)$.
Moreover, suppose that $\Gamma$ contains the pseudo chart
as shown in Fig.~\ref{fig26}$($a$)$,
where $e_1',e_1'',e_2',e_2'',e_3',e_3''$
are internal edges $($possibly terminal edges$)$
of label $m+1$ at $w_1,w_1,w_2,w_2,w_3,w_3$,
respectively.
Then we have the following:

\begin{enumerate}
\item[{\rm (a)}]
$e_1'\not=e_2''$, $e_1''\not=e_2'$ 
$($see Fig.~\ref{fig26}$($b$))$.
\item[{\rm (b)}]
$e_1'\not=e_3''$, $e_1''\not=e_3'$
$($see Fig.~\ref{fig26}$($c$))$.
\end{enumerate}
\end{lemma}

\begin{Proof}
Let $e$ be the terminal edge of label $m$ at $w_1$.
Let $D$ be the special $2$-angled disk of $\Gamma_m$
with $\partial D\ni w_1,w_2$.
Let $e'$ be the internal edge of label $m$
with $\partial e'=\{w_2,w_3\}$.

Without loss of generality,
we can assume that 
\begin{enumerate}
\item[(1)] the terminal edge $e$ is oriented inward at $w_1$.
\end{enumerate}
Then by Assumption~\ref{AssumeTerminal},
\begin{enumerate}
\item[(2)] the two internal edges in $\partial D$
are oriented from $w_2$ to $w_3$,
\item[(3)] the edge $e'$ is oriented from $w_2$ to $w_3$.
\end{enumerate}

We shall show Statement~(a).
Suppose $e_1'=e_2''$
(see Fig.~\ref{fig26}(b)).
Then the edge $e_1'$ separates $Cl(S^2-D)$ into two disks.
One of the two disks contains the terminal edge $e$,
say $F$.
Since $\Gamma$ is of type $(m;5,2)$,
the interior ${\rm Int}F$ does not contain
any white vertices in $\cup^{\infty}_{i=0}\Gamma_{m-i}$.
Thus by Corollary~\ref{CorDiskLemmaWhiteVertice},
we can assume that
$\Gamma_{m-1}\cap e_1'=\emptyset$.
Since the two internal edges in $\partial D$
are oriented inward at $w_2$ by (2),
we can apply Lemma~\ref{C-I-M2A11B12}(a) for the edge $e_1'$ by
(1).
Hence the chart $\Gamma$ is not minimal.
However, this contradicts the fact that 
$\Gamma$ is minimal.
Hence $e_1'\not=e_2''$.

Similarly, we can show $e_1''\not=e_2'$.
Therefore, Statement~(a) holds.

We shall show Statement~(b).
Suppose $e_1'=e_3''$
(see Fig.~\ref{fig26}(c)).
By the similar way as above,
we can assume that
$\Gamma_{m-1}\cap e_1'=\emptyset$.
Since the edge $e'$ is oriented inward at $w_3$ by (3),
we can apply Lemma~\ref{C-I-M2A11B12}(a) for the edge $e_1'$ by
(1).
Hence the chart $\Gamma$ is not minimal.
However, this contradicts the fact that 
$\Gamma$ is minimal.
Hence $e_1'\not=e_3''$.

Similarly, we can show $e_1''\not=e_3'$.
Therefore, Statement~(b) holds.
\end{Proof}

\begin{figure}[htb]
\centerline{\includegraphics{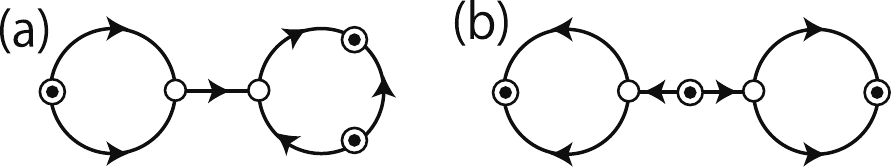}}
\caption{\label{fig27}
Connected components of $\Gamma_m$ with five white vertices.}
\end{figure}

\begin{lemma}
\label{OriGammaM5}
{\rm (\cite[Lemma 7.2(c)]{ChartAppVIII})}
Let $\Gamma$ be a minimal chart, and $m$ a label of $\Gamma$.
Let $G$ be a connected component of $\Gamma_m$
with $w(G)=5$.
If $G$ is the graph as shown in Fig.~\ref{fig10}$($g$)$
$($resp. Fig.~\ref{fig10}$($h$))$,
then $G$ is one of the RO-family of the graph as shown
in Fig.~\ref{fig27}$($a$)$
$($resp. Fig.~\ref{fig27}$($b$))$.
\end{lemma}

\begin{figure}[htb]
\centerline{\includegraphics{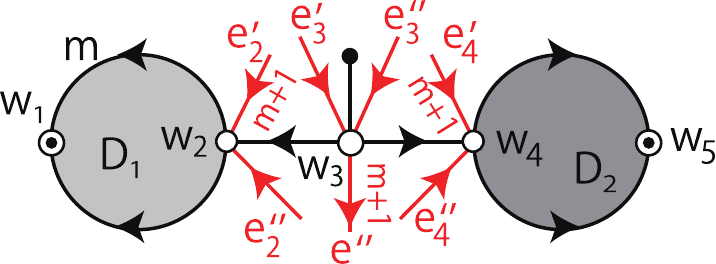}}
\caption{\label{fig28}
The graph as shown in Fig.~\ref{fig10}(h).}
\end{figure}

Let $\Gamma$ be a minimal chart of type $(m;5,2)$.
Suppose that $\Gamma_m$ contains the graph $G$ as
shown in Fig.~\ref{fig10}(h).
Then by Lemma~\ref{OriGammaM5},
the graph $G$ is the graph as shown in 
Fig.~\ref{fig27}(b).
Thus the chart $\Gamma$ contains the pseudo chart 
as shown in Fig.~\ref{fig28}.

From now on throughout this section,
we use the notations as shown in Fig.~\ref{fig28},
where 
\begin{enumerate}
\item[(a)] $w_1,w_2,\cdots,w_5$ are the five white vertices 
in $G$,
\item[(b)] $e_2',e_2'',e_3',e_3'',e_4',e_4''$ are
 internal edges (possibly terminal edges) of label $m+1$
oriented inward at $w_2,w_2,w_3,w_3,w_4,w_4$,
respectively.
\end{enumerate}
Moreover, none of $e_2',e_2'',e_3',e_3'',e_4',e_4''$ are
middle at $w_2,w_3$ or $w_4$,
by Assumption~\ref{AssumeTerminal}
\begin{enumerate}
\item[(c)] none of the six edges 
$e_2',e_2'',e_3',e_3'',e_4',e_4''$ are terminal edges.
\end{enumerate}

\begin{lemma}
\label{TypeHD1D2NoFeeler}
Let $\Gamma$ be a minimal chart of type $(m;5,2)$.
Suppose that $\Gamma_m$ contains the graph $G$ as shown in 
Fig.~\ref{fig10}$($h$)$.
Let $D_1,D_2$ be the special $2$-angled disks of $\Gamma_m$.
Then 
\begin{enumerate}
\item[{\rm (a)}] $w(\Gamma\cap(S^2-(G\cup D_1\cup D_2)))\geqq1$, 
and
\item[{\rm (b)}]
neither $D_1$ nor $D_2$ has a feeler.
\end{enumerate}
\end{lemma}

\begin{figure}[htb]
\centerline{\includegraphics{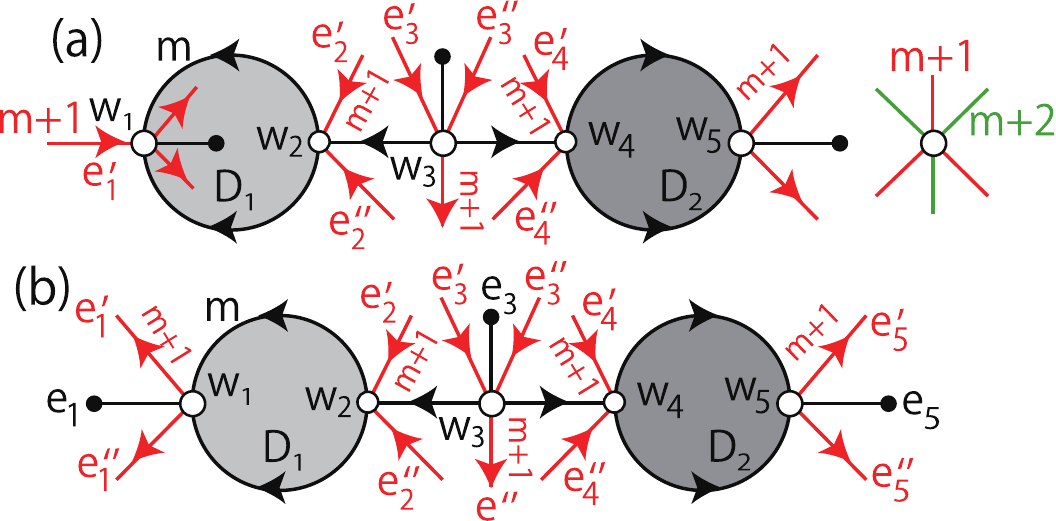}}
\caption{\label{fig29}
(a) The disk $D_1$ has one feeler.
(b) Neither $D_1$ nor $D_2$ has a feeler.}
\end{figure}

\begin{Proof}
We use the notations as shown in Fig.~\ref{fig28}.
By Conditions~(b), (c) of this section,
we have $w(\Gamma\cap(S^2-(G\cup D_1\cup D_2)))\geqq1$
by IO-Calculation with respect to $\Gamma_{m+1}$
in $Cl(S^2-(D_1\cup D_2))$.
Thus Statement~(a) holds.

We shall show Statement~(b).
Suppose that one of $D_1$ and $D_2$ has a feeler.
Without loss of generality we can assume that
$D_1$ has a feeler.
Hence, by Lemma~\ref{Theorem2AngledDisk}
we have $w(\Gamma\cap{\rm Int}D_1)\geqq1$.
Thus by Lemma~\ref{TypeHD1D2NoFeeler}(a),
the condition $w(\Gamma)=7$ implies that
\begin{enumerate}
\item[(1)] $w(\Gamma\cap(S^2-(G\cup D_1\cup D_2)))=1$
and $w(\Gamma\cap{\rm Int}D_2)=0$.
\end{enumerate}
Hence by Lemma~\ref{Theorem2AngledDisk},
the disk $D_2$ has no feeler.
Thus the chart $\Gamma$ contains
the pseudo chart as shown in Fig.~\ref{fig29}(a),
where
\begin{enumerate}
\item[(2)]
$e_1'$ is an internal edge
(possibly a terminal edge) of label $m+1$
oriented inward at $w_1$.
\end{enumerate}
Hence by Condition~(b) of this section,
the seven edges $e_1',e_2',e_2'',e_3',e_3'',e_4',e_4''$
are oriented inward at $w_1,w_2,w_2,w_3,w_3,w_4,w_4$,
respectively.
Thus by Condition~(c) of this section
and by IO-Calculation with respect to $\Gamma_{m+2}$
in $Cl(S^2-(D_1\cup D_2))$,
we have $w(\Gamma_{m+2}\cap(S^2-(G\cup D_1\cup D_2)))\geqq2$.
This contradicts~(1).
Therefore neither $D_1$ nor $D_2$ has a feeler.
Thus Statement~(b) holds.
\end{Proof}

From now on throughout this section,
we use the notations as shown in Fig.~\ref{fig28} 
and Fig.~\ref{fig29}(b), where

\begin{enumerate}
\item[(d)] $e_1',e_1'',e_5',e_5''$ are
 internal edges (possibly terminal edges) of label $m+1$
oriented outward at $w_1,w_1,w_5,w_5$,
respectively,
\item[(e)] $e_1,e_3,e_5$ are terminal edges of label $m$
at $w_1,w_3,w_5$, respectively.
\end{enumerate}


\begin{lemma}
\label{TypeHE1'E5''}
Let $\Gamma$ be a minimal chart of type $(m;5,2)$.
Suppose that $\Gamma_m$ contains the graph as shown 
in Fig.~\ref{fig10}$($h$)$.
If $\Gamma$ contains the pseudo chart as shown in 
Fig.~\ref{fig29}$($b$)$,
then we have the following:

\begin{enumerate}
\item[{\rm (a)}]
$e_3'\not=e_1'$, $e_3''\not=e_5'$,
\item[{\rm (b)}]
$e_3'\not=e_5'$, $e_3''\not=e_1'$, 
\item[{\rm (c)}]
$e_3'\not=e_5''$, $e_3''\not=e_1''$,
\item[{\rm (d)}]
$e_2'\not=e_5'$, $e_4'\not=e_1'$,
\item[{\rm (e)}]
$e_2'\not=e_5''$, $e_4'\not=e_1''$.
\item[{\rm (f)}]
$e_2''\not=e_5'$, $e_2''\not=e_5''$.
\end{enumerate}
\end{lemma}

\begin{figure}[htb]
\centerline{\includegraphics{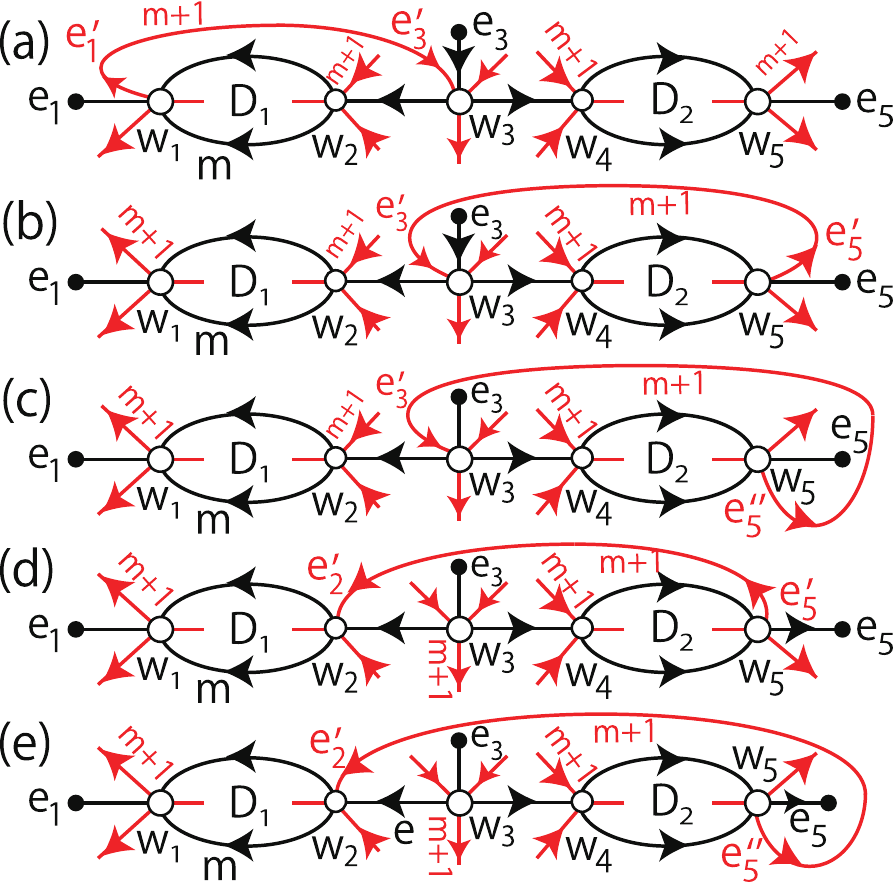}}
\caption{\label{fig30}
(a) $e_3'=e_1'$,
(b) $e_3'=e_5'$,
(c) $e_3'=e_5''$,
(d) $e_2'=e_5'$,
(e) $e_2'=e_5''$.}
\end{figure}

\begin{Proof}
We shall show Statement~(a).
Suppose $e_3'=e_1'$ (see Fig.~\ref{fig30}(a)).
By the similar way of the proof of Lemma~\ref{LemmaTypeGH},
we can assume  that $\Gamma_{m-1}\cap e_3'=\emptyset$.
Since the two edges $e_1,e_3$ are terminal edges at $w_1,w_3$,
respectively,
we can apply Lemma~\ref{C-I-M2A11B12}(d) for the edge $e_3'$.
Thus the chart $\Gamma$ is not minimal.
This contradicts the fact that
$\Gamma$ is minimal.
Hence $e_3'\not=e_1'$.

Similarly, we can show $e_3''\not=e_5'$.
Thus Statement~(a) holds.

We shall show Statement~(b).
Suppose $e_3'=e_5'$ (see Fig.~\ref{fig30}(b)).
By the similar way of the proof of Lemma~\ref{LemmaTypeGH},
we can assume  that $\Gamma_{m-1}\cap e_3'=\emptyset$.
Since the terminal edge $e_3$ is oriented inward at $w_3$,
and since the two internal edges in $\partial D_2$
are oriented inward at $w_5$,
we can apply Lemma~\ref{C-I-M2A11B12}(a) for the edge $e_3'$.
Thus the chart $\Gamma$ is not minimal.
This contradicts the fact that
$\Gamma$ is minimal.
Thus $e_3'\not=e_5'$.

Similarly, we can show $e_3''\not=e_1'$.
Thus Statement~(b) holds.

By the similar way of the proof of Statement~(a), 
we can show Statement~(c) 
(see Fig.~\ref{fig30}(c)).

By the similar way of the proof of Statement~(b), 
we can show Statement~(d)
(see Fig.~\ref{fig30}(d)).

We shall show Statement~(e).
Let $e$ be the internal edge of label $m$ 
with $\partial e=\{w_2,w_3\}$.
Then 
\begin{enumerate}
\item[(1)]
the edge $e$ is not middle at $w_3$.
\end{enumerate}

Suppose $e_2'=e_5''$ (see Fig.~\ref{fig30}(e)).
By the similar way of the proof of Lemma~\ref{LemmaTypeGH},
we can assume  that $\Gamma_{m-1}\cap e_2'=\emptyset$.
Since there exists a terminal edge $e_5$ of label $m$ at $w_5$,
we can apply Lemma~\ref{C-I-M2A11B12}(b) for the edge $e_3'$
by (1).
Thus the chart $\Gamma$ is not minimal.
This contradicts the fact that
$\Gamma$ is minimal.
Thus $e_2'\not=e_5''$.

Similarly, we can show $e_4'\not=e_1''$.
Thus Statement~(e) holds.

By the similar way of the proofs of 
Lemma~\ref{TypeHE1'E5''}(d),(e),
we can show Statement~(f).
\end{Proof}

\begin{lemma}
\label{TypeHE1'E1''WhiteVertex}
Let $\Gamma$ be a minimal chart of type $(m;5,2)$.
If $\Gamma$ contains the pseudo chart as shown in 
Fig.~\ref{fig29}$($b$)$,
then each of $e_3',e_3''$ contains a white vertex
different from the five white vertices $w_1,w_2,\cdots,w_5$.
\end{lemma}

\begin{Proof}
Since $e_3'$ is not a terminal edge by
Condition~(c),
there are six cases:
(i) $e_3'$ is a loop,
(ii) $e_3'=e_1'$,
(iii) $e_3'=e_1''$,
(iv) $e_3'=e_5'$,
(v) $e_3'=e_5''$,
(vi) $e_3'$ contains a white vertex
different from the five white vertices $w_1,w_2,\cdots,w_5$.

By Lemma~\ref{LemmaNoLoop}, 
Case (i) does not occur.
By Lemma~\ref{TypeHE1'E5''}(a),
Case (ii) does not occur.
By Lemma~\ref{LemmaTypeGH}(b),
Case (iii) does not occur.
By Lemma~\ref{TypeHE1'E5''}(b),
Case (iv) dose not occur.
By Lemma~\ref{TypeHE1'E5''}(c),
Case (v) dose not occur.
Therefore, Case (vi) occurs.

Similarly, we can show that $e_3''$ contains a white vertex
different from $w_1,w_2,\cdots,w_5$.
\end{Proof}


\begin{lemma}
\label{TypeHE2'E1}
Let $\Gamma$ be a minimal chart of type $(m;5,2)$.
Suppose that $\Gamma_m$ contains the graph as shown 
in Fig.~\ref{fig10}$($h$)$.
If $\Gamma$ contains the pseudo chart as shown in 
Fig.~\ref{fig29}$($b$)$,
then $e_2'\not\ni w_3$
 and
$e_4'\not\ni w_3$.
\end{lemma}

\begin{figure}[htb]
\centerline{\includegraphics{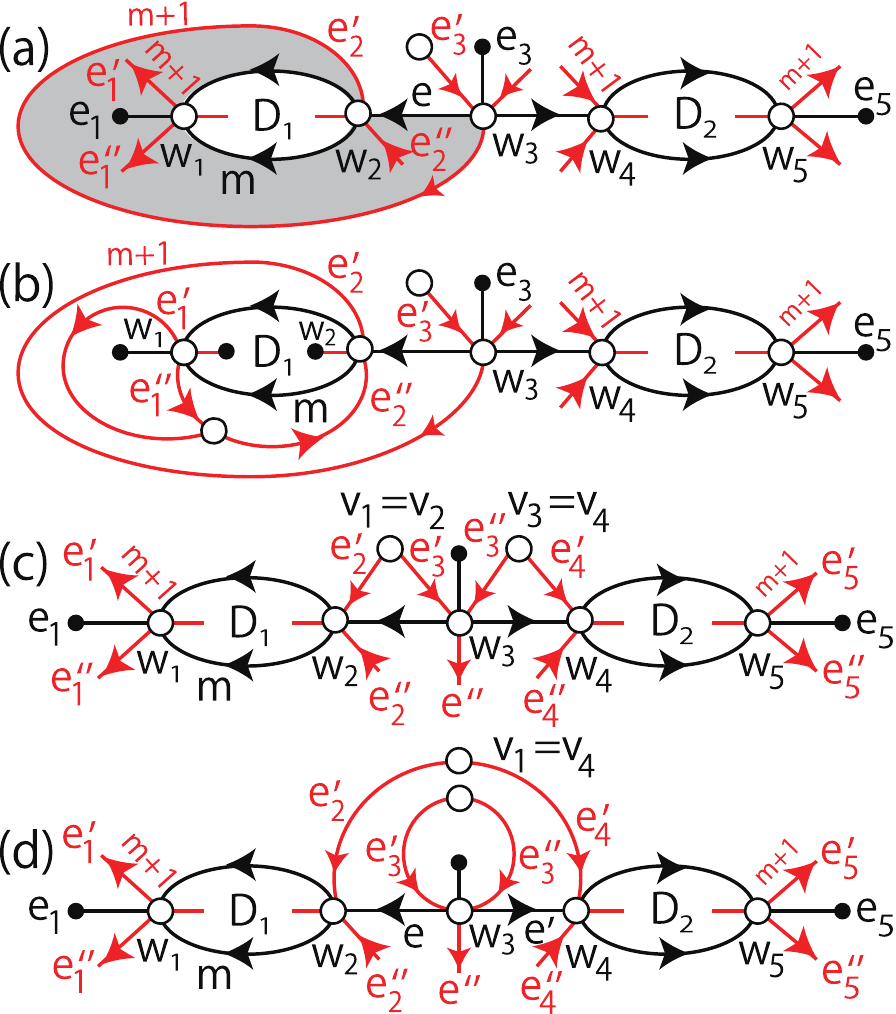}}
\caption{\label{fig31}
(a) The gray region is the region $F_1$.
(b) $e_1'\cap e_1''\cap e_2''$ is a white vertex.
(c) $v_1=v_2$, $v_3=v_4$.
(d) $v_1=v_4$, $v_2=v_3$. }
\end{figure}

\begin{Proof}
Let $e$ be the internal edge of label $m$
with $\partial e=\{w_2,w_3\}$.

Suppose $e_2'\ni w_3$ (see Fig.~\ref{fig31}(a)).
Then the curve $e_2'\cup e$
separates $Cl(S^2-(D_1\cup D_2))$ into two regions.
One of the two regions contains the edge $e_2''$,
say $F_1$.
Let $F_2$ be the other region.

By Lemma~\ref{TypeHE1'E1''WhiteVertex},
the edge $e_3'$ contains a white vertex different from
the five white vertices $w_1,w_2,\cdots,w_5$.
Thus  $w(\Gamma\cap{\rm Int}F_2)\geqq1$.

Next, we shall show that 
the edge $e_1''$ contains a white vertex in ${\rm Int}F_1$.
Since the edge $e_1''$ is not middle at $w_1$,
by Assumption~\ref{AssumeTerminal}
the edge $e_1''$ is not a terminal edge.
Since the edges $e_1',e_1''$ is oriented outward at $w_1$,
we have $e_1''\not=e_1'$.
Hence either $e_1''=e_2''$ or 
$e_1''$ contains a white vertex in ${\rm Int}F_1$.
If $e_1''=e_2''$,
then there exists a lens.
This contradicts Lemma~\ref{NoLens}.
Thus the edge $e_1''$ contains a white vertex in ${\rm Int}F_1$.

Let $w_6$ be the white vertex in ${\rm Int}F_1$ 
with $w_6\in e_1''$.
Since $w(\Gamma\cap{\rm Int}F_2)\geqq1$ and $w(\Gamma)=7$,
\begin{enumerate}
\item[(1)] 
${\rm Int}F_1$ contains exactly one white vertex $w_6$,
\item[(2)] $w(\Gamma\cap{\rm Int}D_1)=0$.
\end{enumerate}

Next, we shall show $e_1'\ni w_6$.
Similarly, we can show that 
the edge $e_1'$ is not a terminal edge.
Hence either $e_1'=e_2''$ or $e_1''\ni w_6$.
If $e_1'=e_2''$, then 
this contradicts Lemma~\ref{LemmaTypeGH}(a).
Thus $e_1'\ni w_6$.

By Condition~(c) of this section,
the edge $e_2''$ is not a terminal edge.
Hence we have $w_6\in e_1'\cap e_1''\cap e_2''$.
Moreover, by (2) and Lemma~\ref{Theorem2AngledDisk},
a regular neighborhood of $D_1$ contains
the pseudo chart as shown in Fig.~\ref{fig03}(b).
Thus the chart $\Gamma$ contains the pseudo chart as
shown in Fig.~\ref{fig24}
(see Fig.~\ref{fig31}(b)).
Hence by Lemma~\ref{CorDiskLemma},
we have $w(\Gamma\cap{\rm Int}F_1)\geqq2$.
This contradicts (1).
Thus $e_2'\not \ni w_3$.

Similarly we can show $e_4'\not\ni w_3$.
\end{Proof}

\begin{lemma}
\label{TypeHE2'E3'WhiteVertex}
Let $\Gamma$ be a minimal chart of type $(m;5,2)$.
If $\Gamma$ contains the pseudo chart as shown in 
Fig.~\ref{fig29}$($b$)$,
then each of $e_2',e_4'$ contains a white vertex
different from the five white vertices $w_1,w_2,\cdots,w_5$.
\end{lemma}

\begin{Proof}
Since $e_2'$ is not a terminal edge by
Condition~(c),
there are six cases:
(i) $e_2'\ni w_3$,
(ii) $e_2'=e_1'$,
(iii) $e_2'=e_1''$,
(iv) $e_2'=e_5'$,
(v) $e_2'=e_5''$,
(vi) $e_2'$ contains a white vertex
different from the five white vertices $w_1,w_2,\cdots,w_5$.

By Lemma~\ref{TypeHE2'E1},
Case (i) does not occur.
By Lemma~\ref{NoLens},
Case (ii) does not occur.
By Lemma~\ref{LemmaTypeGH}(a),
Case (iii) does not occur.
By Lemma~\ref{TypeHE1'E5''}(d),
Case (iv) does not occur.
By Lemma~\ref{TypeHE1'E5''}(e),
Case (v) does not occur.
Therefore Case (vi) occurs.

Similarly, we can show that $e_4'$ contains a white vertex
different from $w_1,w_2,\cdots,w_5$.
\end{Proof}

\begin{proposition}
\label{NoTypeH}
Let $\Gamma$ be a minimal chart of type $(m;5,2)$.
Then $\Gamma_m$ does not contain the graph as shown
in Fig.~\ref{fig10}$($h$)$.
\end{proposition}

\begin{Proof}
Suppose that $\Gamma_m$ contains the graph $G$ as shown
in Fig.~\ref{fig10}(h).
Then by Lemma~\ref{TypeHD1D2NoFeeler}(b),
we can assume that
the chart $\Gamma$ contains
the pseudo chart as shown in Fig.~\ref{fig29}(b).
We use the notations as shown in Fig.~\ref{fig29}(b),
where 
\begin{enumerate}
\item[(1)] the five edges $e'',e_1',e_1'',e_5',e_5''$
are internal edges (possibly terminal edges) of label $m+1$ 
oriented outward at $w_3,w_1,w_1,w_5,w_5$,
respectively.
\end{enumerate}
Moreover, none of $e_1',e_1'',e_5',e_5''$
are middle at $w_1$ or $w_5$.
Thus by Assumption~\ref{AssumeTerminal},
\begin{enumerate}
\item[(2)] none of $e_1',e_1'',e_5',e_5''$ are terminal edges.
\end{enumerate}

By Lemma~\ref{TypeHE1'E1''WhiteVertex} and 
Lemma~\ref{TypeHE2'E3'WhiteVertex},
each of $e_2',e_3',e_3'',e_4'$
contains a white vertex different from
the five white vertices $w_1,w_2,\cdots,w_5$.
Let $v_1,v_2,v_3,v_4$ be white vertices
different from $w_1,w_2,\cdots,w_5$ with
$v_1\in e_2',v_2\in e_3',v_3\in e_3'',v_4\in e_4'$.
Then the condition $w(\Gamma)=7$ implies
that 
\begin{enumerate}
\item[(3)] $w(\Gamma\cap(S^2-(G\cup D_1\cup D_2)))=2$.
\end{enumerate}
Hence, the set $\{v_1,v_2,v_3,v_4\}$ consists
of two white vertices.
Since by Condition~(b) the four edges  $e_2',e_3',e_3'',e_4'$
are oriented inward at $w_2,w_3,w_3,w_4$,
respectively,
there are two cases:
(i) $v_1=v_2$, $v_3=v_4$  
(see Fig.~\ref{fig31}(c)),
(ii) $v_1=v_4$, $v_2=v_3$
(see Fig.~\ref{fig31}(d)).

{\bf Case (i).}
By (1), (2), (3),
\begin{enumerate}
\item[(4)] the edge $e''$ must be a terminal edge.
\end{enumerate}
Moreover,
since the edge $e_2''$ is not a terminal edge
by Condition~(c) in this section,
there are five cases:
(i-1) $e_2''=e_1'$,
(i-2) $e_2''=e_1''$,
(i-3) $e_2''=e_5'$,
(i-4) $e_2''=e_5''$,
(i-5) $e_2''\ni v_1$ or $e_2''\ni v_3$.

By Lemma~\ref{LemmaTypeGH}(a), Case (i-1) does not occur.
By Lemma~\ref{NoLens}, Case (i-2) does not occur.
By Lemma~\ref{TypeHE1'E5''}(f),
neither Case (i-3) nor Case (i-4) occurs.

For Case (i-5),
if $e_2''\ni v_1$,
then there exist three internal edges $e_2',e_2'',e_3'$
of label $m+1$ oriented outward at $v_1$.
This contradicts the definition of the chart.
Similarly, 
if $e_2''\ni v_3$,
then we have the same contradiction.
Thus Case (i-5) does not occur.

Hence all the five cases 
do not occur. Thus Case (i) does not occur.

{\bf Case (ii).}
Let $e,e'$ be the internal edges of label $m$
with $\partial e=\{w_2,w_3\}$,
$\partial e'=\{w_3,w_4\}$.
Since $v_1=v_4$,
the curve $e_2'\cup e_4'\cup e\cup e'$
separates $Cl(S^2-(D_1\cup D_2))$ into two regions.
One of the two regions contains the edge $e_1'$,
say $F$.
By (1), (2) and IO-Calculation
with respect to $\Gamma_{m+1}$ in $F$,
we have $w(\Gamma\cap{\rm Int }F)\geqq 1$.
Thus $w(\Gamma\cap(S^2-(G\cup D_1\cup D_2)))\geqq3$.
This contradicts (3).
Hence Case (ii) does not occur.

Therefore both cases (i),(ii) do not occur.
Hence $\Gamma_m$ does not contain the graph as shown
in Fig.~\ref{fig10}(h).
\end{Proof}



\section{Case of the graph as shown in Fig.~\ref{fig10}(g)}
\label{s:TypeG}

In this section,
we shall show the main theorem.

\begin{figure}[htb]
\centerline{\includegraphics{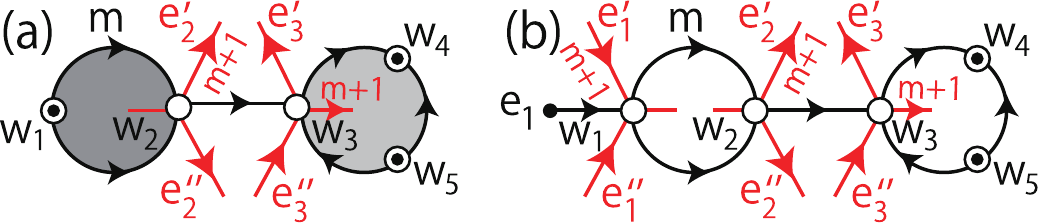}}
\caption{\label{fig32}
The graphs as shown in Fig.~\ref{fig10}(g).
(a) The light gray region is the disk $D_1$.
The dark gray region is the disk $D_2$.
(b) The disk $D_2$ has no feeler.}
\end{figure}

Let $\Gamma$ be a minimal chart of type $(m;5,2)$.
Suppose that $\Gamma_m$ contains the graph $G$ as
shown in Fig.~\ref{fig10}(g).
Then by Lemma~\ref{OriGammaM5},
the graph $G$ is the graph as shown in 
Fig.~\ref{fig27}(a).
Thus the chart $\Gamma$ contains the pseudo chart 
as shown in Fig.~\ref{fig32}(a).

From now on throughout this section,
we use the notations as shown in Fig.~\ref{fig32}(a),
where 
\begin{enumerate}
\item[(a)] $w_1,w_2,\cdots,w_5$ are the five white vertices 
in $G$,
\item[(b)] $D_1$ is a special 3-angled disk of $\Gamma_m$,
and $D_2$ is a special 2-angled disk of $\Gamma_m$,
\item[(c)] $e_2',e_2'',e_3'$
are internal edges (possibly terminal edges) of label $m+1$ 
oriented outward at $w_2,w_2,w_3$, respectively.
\end{enumerate}
In particular, if $D_2$ has no feeler,
then the chart $\Gamma$ contains the pseudo chart 
as shown in Fig.~\ref{fig32}(b),
where 
\begin{enumerate}
\item[(d)] $e_1',e_1''$
are internal edges (possibly terminal edges) of label $m+1$ 
oriented inward at $w_1$.
\end{enumerate}
Moreover, none of $e_1',e_1'',e_2',e_2'',e_3'$ 
are middle at $w_1,w_2$ or $w_3$.
Thus by Assumption~\ref{AssumeTerminal},
\begin{enumerate}
\item[(e)] none of $e_1',e_1'',e_2',e_2'',e_3'$ 
are terminal edges.
\end{enumerate}

\begin{lemma}
\label{TypeGE1E1}
Let $\Gamma$ be a minimal chart of type $(m;5,2)$.
Suppose that $\Gamma_m$ contains the graph $G$ as shown in
Fig.~\ref{fig10}$($g$)$.
If the $2$-angled disk $D_2$ has no feeler,
then $e_1'\not=e_3'$.
\end{lemma}

\begin{figure}[htb]
\centerline{\includegraphics{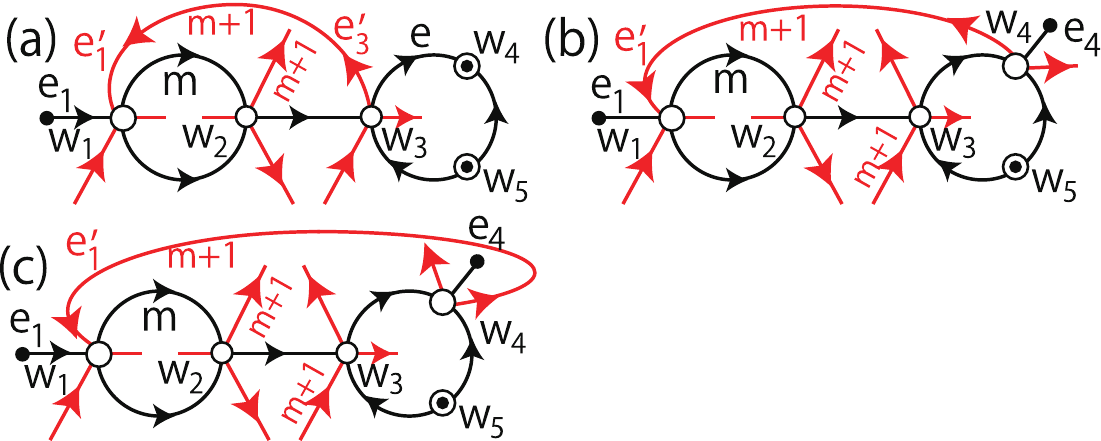}}
\caption{\label{fig33}
(a) $e_1'=e_3'$, (b),(c) $e_1'\ni w_4$.}
\end{figure}

\begin{Proof}
Let $e$ be the internal edge of label $m$ with
$\partial e=\{w_3,w_4\}$.
Then 
\begin{enumerate}
\item[(1)] the edge $e$ is not middle at $w_4$.
\end{enumerate}

Suppose $e_1'=e_3'$ (see Fig.~\ref{fig33}(a)).
By the similar way of the proof of Lemma~\ref{LemmaTypeGH}(a),
we can assume that $\Gamma_{m-1}\cap e_1'=\emptyset$.
Since there exists a terminal edge of label $m$ at $w_1$,
we can apply Lemma~\ref{C-I-M2A11B12}(b) 
for the edge $e_1'$ by (1). 
Hence the chart $\Gamma$ is not minimal.
However, this contradicts the fact that 
$\Gamma$ is minimal.
Thus $e_1'\not=e_3'$.
\end{Proof}

\begin{lemma}
\label{TypeGE1W4}
Let $\Gamma$ be a minimal chart of type $(m;5,2)$.
Suppose that $\Gamma_m$ contains the graph $G$ as shown in
Fig.~\ref{fig10}$($g$)$.
If the $2$-angled disk $D_2$ has no feeler,
then $e_1'\not \ni w_4$ and
$e_1''\not \ni w_4$.
\end{lemma}

\begin{Proof}
Let $e_4$ be the terminal edge of label $m$ at $w_4$,
and $D_1$ the special 3-angled disk of $\Gamma_m$
with $\partial D_1\subset G$.
Then 
\begin{enumerate}
\item[(1)] the edge $e_4$ is oriented outward at $w_4$
(see Fig.~\ref{fig32}(b)).
\end{enumerate}

We shall show  $e_1'\not\ni w_4$.
Suppose $e_1'\ni w_4$.
Then by (1),
we have $e_4\not\subset D_1$
(see Fig.~\ref{fig33}(b),(c)).
By the similar way of the proof of Lemma~\ref{LemmaTypeGH}(a),
we can assume that $\Gamma_{m-1}\cap e_1'=\emptyset$.

If the chart $\Gamma$ contains the pseudo chart
as shown in Fig.~\ref{fig33}(b),
then we can apply Lemma~\ref{C-I-M2A11B12}(d) 
for the edge $e_1'$.
Thus $\Gamma$ is not minimal.
However, this contradicts the fact that 
$\Gamma$ is minimal.
Hence the chart $\Gamma$ does not contain the pseudo chart
as shown in Fig.~\ref{fig33}(b).

If the chart $\Gamma$ contains the pseudo chart
as shown in Fig.~\ref{fig33}(c),
then we have the same contradiction
by the similar way of the proof of Lemma~\ref{LemmaTypeGH}(a).
Thus the chart $\Gamma$ does not contain the pseudo chart
as shown in Fig.~\ref{fig33}(c).
Therefore $e_1'\not\ni w_4$.

Similarly we can show $e_1''\not\ni w_4$.
\end{Proof}


By the similar way of the proof of Lemma~\ref{TypeGE1W4},
we can show the following lemma.

\begin{lemma}
\label{TypeGE5E5}
Let $\Gamma$ be a minimal chart of type $(m;5,2)$.
Suppose that $\Gamma_m$ contains the graph $G$ as shown in
Fig.~\ref{fig10}$($g$)$.
Let $D_1$ be the special $3$-angled disk of $\Gamma_m$
with $\partial D_1\subset G$.
Let $e_5$ be the terminal edge of label $m$ 
 oriented inward at $w_5$,
and $e_5',e_5''$ internal edges $($possibly terminal edges$)$
of label $m+1$ oriented inward at $w_5$.
If $e_5\not\subset D_1$,
then we have the following:

\begin{enumerate}
\item[{\rm (a)}] $e_5'\not\ni w_2,e_5''\not\ni w_2$,
\item[{\rm (b)}] $e_5'\not\ni w_3,e_5''\not\ni w_3$,
\item[{\rm (c)}] $e_5'\not\ni w_4,e_5''\not\ni w_4$.
\end{enumerate}
\end{lemma}


\begin{lemma}
\label{TypeGTwoFeeler}
Let $\Gamma$ be a minimal chart of type $(m;5,2)$.
Suppose that $\Gamma_m$ contains the graph $G$ as shown in
Fig.~\ref{fig10}$($g$)$.
Let $D_1$ be the special $3$-angled disk of $\Gamma_m$
with $\partial D_1\subset G$.
Then $D_1$ contains at most one feeler.
\end{lemma}

\begin{figure}[htb]
\centerline{\includegraphics{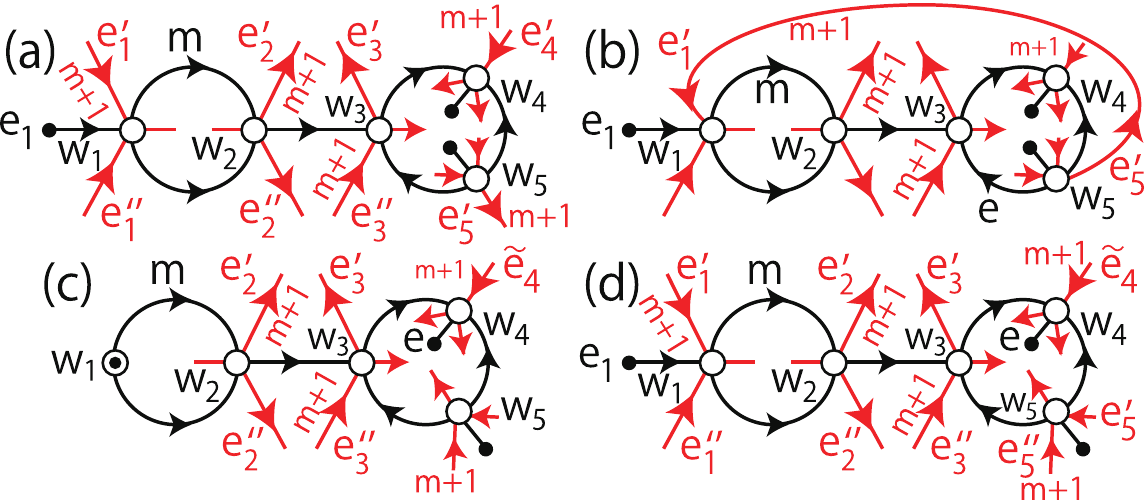}}
\caption{\label{fig34}
(a) The $3$-angled disk $D_1$ has two feelers.
(b) $e_1'\ni w_5$.
(c) The $3$-angled disk $D_1$ has 
exactly one feeler $e$ at $w_4$.
(d) The 2-angled disk $D_2$ has no feeler.}
\end{figure}

\begin{Proof}
Let $D_2$ be the special 2-angled disk of $\Gamma_m$
with $\partial D_2\subset G$.

Suppose that the 3-angled disk $D_1$ 
contains at least two feelers.
Then $D_1$ contains exactly two feelers.
Thus by Lemma~\ref{Theorem3AngledDiskChart9},
we have $w(\Gamma\cap{\rm Int}D_1)\geqq 2$.
Hence the condition $w(\Gamma)=7$ implies that
\begin{enumerate}
\item[(1)] $w(\Gamma\cap(S^2-(D_1\cup D_2)))=0$, 
$w(\Gamma\cap{\rm Int}D_2)=0$.
\end{enumerate}
Thus by Lemma~\ref{Theorem2AngledDisk},
the disk $D_2$ has no feeler.
Hence the chart $\Gamma$ contains the pseudo chart as shown
in Fig.~\ref{fig34}(a).
We use the notations as shown in Fig.~\ref{fig34}(a),
where 
$e_1'$ is an internal edge (possibly a terminal edge)
of label $m+1$ oriented inward at $w_1$.
By Condition~(e) of this section,
the edge $e_1'$ is not a terminal edge.
Thus there are four cases:
(i) $e_1'=e_2'$,
(ii) $e_1'=e_2''$,
(iii) $e_1'=e_3'$,
(iv) $e_1'=e_5'$ (see Fig.~\ref{fig34}(b)).

By Lemma~\ref{NoLens},
Case (i) does not occur.
By Lemma~\ref{LemmaTypeGH}(a),
Case (ii) does not occur.
By Lemma~\ref{TypeGE1E1},
Case (iii) does not occur.
Hence we shall show that Case (iv) does not occur.

{\bf Case (iv).}
Let $e$ be the internal edge of label $m$
with $\partial e=\{w_3,w_5\}$.
Then 
\begin{enumerate}
\item[(2)] the edge $e$ is not middle at $w_3$.
\end{enumerate}
By the similar way of the proof of Lemma~\ref{LemmaTypeGH}(a),
we can assume that $\Gamma_{m-1}\cap e_1'=\emptyset$.
Since there exists a terminal edge $e_1$ of label $m$ at $w_1$,
we can apply Lemma~\ref{C-I-M2A11B12}(b) 
for the edge $e_1'$ by (2). 
Hence the chart $\Gamma$ is not minimal.
However, this contradicts the fact that 
$\Gamma$ is minimal.
Hence Case (iv) does not occur.

Therefore all the four cases do not occur.
Thus the 3-angled disk $D_1$ contains at most one feeler.
\end{Proof}

\begin{lemma}
\label{TypeGOneFeelerW5}
Let $\Gamma$ be a minimal chart of type $(m;5,2)$.
Suppose that $\Gamma_m$ contains the graph $G$ as shown in
Fig.~\ref{fig10}$($g$)$ $($see Fig.~\ref{fig32}$($a$))$.
Let $D_1$ be the special $3$-angled disk of $\Gamma_m$
with $\partial D_1\subset G$.
If $D_1$ contains a feeler $e$,
then the feeler $e$ contains the white vertex $w_5$.
\end{lemma}

\begin{Proof}
Let $D_2$ be the special 2-angled disk of $\Gamma_m$
with $\partial D_2\subset G$.

Suppose that $e\ni w_4$.
Since 
the 3-angled disk $D_1$ has at most one feeler
by Lemma~\ref{TypeGTwoFeeler},
the disk $D_1$ has exactly one feeler $e$.
Thus the chart $\Gamma$ contains the pseudo chart
as shown in Fig.~\ref{fig34}(c).
Let $e_3,e_4',e_4'',\widetilde{e}_5$ be internal edges
(possibly terminal edges) of label $m+1$ oriented
outward at $w_3,w_4,w_4,w_5$, respectively, in $D_1$.
Then none of $e_3,e_4',e_4''$ are middle at $w_3$ or $w_4$.
Thus by Assumption~\ref{AssumeTerminal},
 none of $e_3,e_4',e_4''$ are terminal edges.
Hence by IO-Calculation with respect to $\Gamma_{m+1}$
in $D_1$,
we have $w(\Gamma\cap{\rm Int}D_1)\geqq2$.
Hence, the condition $w(\Gamma)=7$ implies that
\begin{enumerate}
\item[(1)] 
$w(\Gamma\cap(S^2-(D_1\cup D_2)))=0$, 
$w(\Gamma\cap{\rm Int}D_2)=0$.
\end{enumerate}
Thus by Lemma~\ref{Theorem2AngledDisk},
the 2-angled disk $D_2$ has no feeler.
Hence, the chart $\Gamma$ contains the pseudo chart
as shown in Fig.~\ref{fig34}(d), where
$e_1',e_1'',e_3'',\widetilde{e}_4,e_5',e_5''$ are six 
internal edges
(possibly terminal edges) of label $m+1$
oriented inward at $w_1,w_1,w_3,w_4,w_5,w_5$,
respectively.
Moreover,
none of $e_1',e_1'',e_5',e_5''$ are middle at $w_1$ or $w_5$.
Hence, none of $e_1',e_1'',e_5',e_5''$ are terminal edges.
Thus, by IO-Calculation with respect to $\Gamma_{m+1}$
in $Cl(S^2-(D_1\cup D_2))$,
we have $w(\Gamma\cap(S^2-(D_1\cup D_2)))\geqq1$.
This contradicts~(1).
Therefore, the feeler $e$ contains the white vertex $w_5$.
\end{Proof}

\begin{lemma}
\label{TypeGOneFeelerD1OneFeelerD2}
Let $\Gamma$ be a minimal chart of type $(m;5,2)$.
Suppose that $\Gamma_m$ contains the graph $G$ as shown in
Fig.~\ref{fig10}$($g$)$ $($see Fig.~\ref{fig32}$($a$))$.
Let $D_1$ be the special $3$-angled disk of $\Gamma_m$
with $\partial D_1\subset G$.
Then the $3$-angled disk $D_1$ has no feeler.
\end{lemma}

\begin{figure}[htb]
\centerline{\includegraphics{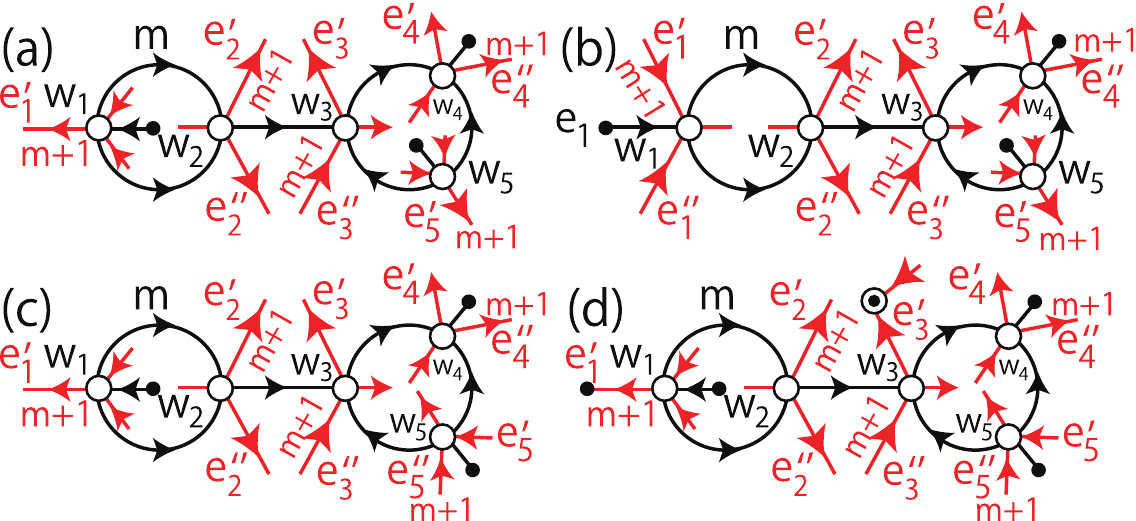}}
\caption{\label{fig35}
(a),(c),(d) The 2-angled disk $D_2$ has one feeler.
(b) The 2-angled disk $D_2$ has no feeler.}
\end{figure}

\begin{Proof}
Suppose that the $3$-angled disk $D_1$ has a feeler $e$.
By Lemma~\ref{Theorem3AngledDisk},
we have 
\begin{enumerate}
\item[(1)] $w(\Gamma\cap{\rm Int}D_1)\geqq1$.
\end{enumerate}
By Lemma~\ref{TypeGTwoFeeler},
the disk $D_1$ has exactly one feeler $e$.
Moreover, by Lemma~\ref{TypeGOneFeelerW5},
the feeler $e$ contains the white vertex $w_5$.

Let $D_2$ be the sepecial $2$-angled disk of $\Gamma_m$
with $\partial D_2\subset G$.
There are two cases:
(i) $D_2$ has one feeler
(see Fig.~\ref{fig35}(a)),
(ii) $D_2$ has no feeler
(see Fig.~\ref{fig35}(b)).

{\bf Case (i).}
We use the notations as shown in Fig.~\ref{fig35}(a),
where
$e_1',e_2',e_2'',e_3',e_4',e_4'',e_5'$ 
are seven internal edges (possibly terminal edges)
of label $m+1$ oriented outward at 
$w_1,w_2,w_2,w_3,w_4,w_4,w_5$, respectively.
Moreover, none of $e_2',e_2'',e_3',e_4',e_4''$
are middle at $w_2,w_3$ or $w_4$.
Thus by Assumption~\ref{AssumeTerminal},
none of  $e_2',e_2'',e_3',e_4',e_4''$ are
terminal edges.
Hence by IO-Calculation with respect to $\Gamma_{m+1}$ in
$Cl(S^2-(D_1\cup D_2))$,
we have  $w(\Gamma\cap(S^2-(D_1\cup D_2)))\geqq2$.
Thus by (1), we have\vspace{2mm}

$\begin{array}{rcl}
7=w(\Gamma) & = & w(G)+w(\Gamma\cap{\rm Int}D_1)+w(\Gamma\cap(S^2-(D_1\cup D_2)))\vspace{2mm}\\
& \geqq & 5+1+2=8.\vspace{2mm}
\end{array}
$\\
This is a contradiction.
Hence Case (i) does not occur.

{\bf Case (ii).}
We use the notations as shown in Fig.~\ref{fig35}(b), where
\begin{enumerate}
\item[(2)] $e_2',e_2'',e_3',e_4',e_4'',e_5'$ 
are six internal edges (possibly terminal edges)
of label $m+1$ oriented outward at 
$w_2,w_2,w_3,w_4,w_4,w_5$, respectively.
\end{enumerate}
Moreover, none of $e_2',e_2'',e_3',e_4',e_4''$
are middle at $w_2,w_3$ or $w_4$.
Thus by Assumption~\ref{AssumeTerminal},
\begin{enumerate}
\item[(3)] 
none of  $e_2',e_2'',e_3',e_4',e_4''$ are
terminal edges.
\end{enumerate}
Hence by IO-Calculation with respect to $\Gamma_{m+1}$ in
$Cl(S^2-(D_1\cup D_2))$,
we have  $w(\Gamma\cap(S^2-(D_1\cup D_2)))\geqq1$.
Thus by (1), the condition $w(\Gamma)=7$ implies
that 
\begin{enumerate}
\item[(4)] $w(\Gamma\cap(S^2-(D_1\cup D_2)))=1$.
\end{enumerate}

Let $w_7$ be the white vertex in $S^2-(D_1\cup D_2)$.
Then by (2),(3),(4),
there are two internal edges of label $m+1$ oriented inward
at $w_7$.
Moreover, there exists a terminal edge of label $m+1$ at $w_7$,
and the edge $e_5'$ must be a terminal edge.

By Condition~(e),
the edge $e_1'$ is not a terminal edge.
Thus there are four cases:
(ii-1) $e_1'=e_2'$ ,
(ii-2) $e_1'=e_2''$,
(ii-3) $e_1'=e_3'$,
(ii-4) $e_1'\ni w_4$ (i.e. $e_1'=e_4'$ or $e_1'=e_4''$).

For Case (ii-1),
there exists a lens.
This contradicts Lemma~\ref{NoLens}.
Thus Case (ii-1) does not occur.
By Lemma~\ref{LemmaTypeGH}(a),
Case (ii-2) does not occur.
By Lemma~\ref{TypeGE1E1},
Case (ii-3) does not occur.
By Lemma~\ref{TypeGE1W4},
Case (ii-4) does not occur.
Therefore, all the four cases do not occur.
Hence Case (ii) does not occur.

Thus both Cases (i),(ii) do not occur.
Therefore $D_1$ has no feeler.
\end{Proof}

\begin{lemma}
\label{TypeGNoFeelerD2}
Let $\Gamma$ be a minimal chart of type $(m;5,2)$.
Suppose that $\Gamma_m$ contains the graph $G$ as shown in
Fig.~\ref{fig10}$($g$)$ $($see Fig.~\ref{fig32}$($a$))$.
Let $D_2$ be the special $2$-angled disk of $\Gamma_m$
with $\partial D_2\subset G$.
Then $D_2$ has no feeler.
\end{lemma}

\begin{Proof}
Let $D_1$ be the special 3-angled disk of $\Gamma_m$
with $\partial D_1\subset G$.

Suppose that $D_2$ has a feeler.
Then $D_2$ has exactly one feeler.
Thus by Lemma~\ref{Theorem2AngledDisk},
we have 
\begin{enumerate}
\item[(1)] $w(\Gamma\cap{\rm Int}D_2)\geqq 1$.
\end{enumerate}
Moreover, 
by Lemma~\ref{TypeGOneFeelerD1OneFeelerD2},
the $3$-angled disk $D_1$ has no feeler.
Thus the chart $\Gamma$ contains the pseudo chart 
as shown in Fig.~\ref{fig35}(c), where
\begin{enumerate}
\item[(2)] $e_1',e_2',e_2'',e_3',e_4',e_4''$ are
internal edges (possibly terminal edges) of label $m+1$
oriented outward at $w_1,w_2,w_2,w_3,w_4,w_4$,
respectively.
\end{enumerate}
Moreover, none of $e_2',e_2'',e_3',e_4',e_4''$ are middle
at $w_2,w_3$ or $w_4$.
Thus by Assumption~\ref{AssumeTerminal},
\begin{enumerate}
\item[(3)] none of $e_2',e_2'',e_3',e_4',e_4''$ are 
terminal edges.
\end{enumerate}
Hence by IO-Calculation with respect to $\Gamma_{m+1}$
in $Cl(S^2-(D_1\cup D_2))$,
we have $w(\Gamma\cap(S^2-(D_1\cup D_2)))\geqq1$.
Thus by (1),
the condition $w(\Gamma)=7$ implies that
\begin{enumerate}
\item[(4)] $w(\Gamma\cap(S^2-(D_1\cup D_2)))=1$.
\end{enumerate}

Let $w_7$ be the white vertex in $S^2-(D_1\cup D_2)$.

{\bf Claim.} The edge $e_3'$ contains the white vertex $w_7$.

{\it Proof of Claim.} 
By (3), the edge $e_3'$ is not a terminal edge.
Hence there are four cases:
$e_3'=e_3''$, $e_3'=e_5'$, $e_3'=e_5''$, or $e_3'\ni w_7$.

If $e_3'=e_3''$,
then the edge $e_3'$ is a loop.
This contradicts Lemma~\ref{LemmaNoLoop}.

If $e_3'=e_5'$ or $e_3'=e_5''$,
then $e_5'\ni w_3$ or $e_5''\ni w_3$.
This contradicts Lemma~\ref{TypeGE5E5}(b).
Hence $e_3'\ni w_7$.
Thus Claim holds.
{\hfill {$\square$}\vspace{1.5em}}

By (2),(3),(4),
the edge $e_1'$ must be a terminal edge of label $m+1$ at $w_1$.
Moreover,
there exists an internal edge of label $m+1$ oriented inward 
at $w_7$ different from $e_3'$, 
and there exists a terminal edge of label $m+1$
at $w_7$ (see Fig.~\ref{fig35}(d)).

Now, the edge $e_5'$ is not middle at $w_5$.
Thus by Assumption~\ref{AssumeTerminal},
the edge $e_5'$ is not a terminal edge.
Hence there are two cases:
$e_5'\ni w_2$ or
$e_5'\ni w_4$.
However this contradicts Lemma~\ref{TypeGE5E5}(a),(c).
Thus the 2-angled disk $D_2$ has no feeler.
\end{Proof}

\begin{proposition}
\label{NoTypeG}
Let $\Gamma$ be a minimal chart of type $(m;5,2)$.
Then $\Gamma_m$ does not contain the graph 
as shown in Fig.~\ref{fig10}$($g$)$.
\end{proposition}

\begin{figure}[htb]
\centerline{\includegraphics{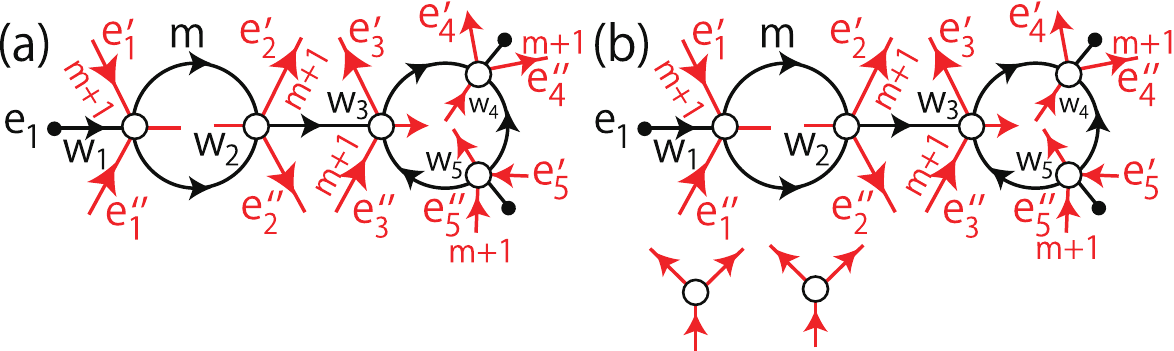}}
\caption{\label{fig36}
 Neither $D_1$ nor $D_2$ has a feeler.}
\end{figure}

\begin{Proof}
Suppose that $\Gamma_m$ contains the graph 
as shown in Fig.~\ref{fig10}$($g$)$,
say $G$.
Let $D_1$ be the special 3-angled disk of $\Gamma_m$,
and $D_2$ the special 2-angled disk of $\Gamma_m$
with $\partial D_1\subset G$ and 
$\partial D_2\subset G$.
By Lemma~\ref{TypeGOneFeelerD1OneFeelerD2} and 
Lemma~\ref{TypeGNoFeelerD2},
neither $D_1$ nor $D_2$ has a feeler.
Moreover, by Lemma~\ref{OriGammaM5},
the graph $\Gamma_m$ contains the graph as shown
in Fig.~\ref{fig27}(a).
Thus the chart $\Gamma$ contains the pseudo chart
as shown in Fig.~\ref{fig36}(a).

{\bf Claim 1.}
The edge $e_1'$ contains a white vertex in $S^2-(D_1\cup D_2)$.

{\it Proof of Claim $1$.}
By Condition~(e) of this section,
the edge $e_1'$ is not a terminal edge.
Thus there are five cases:
$e_1'=e_2'$, $e_1'=e_2''$, $e_1'=e_3'$, $e_1'\ni w_4$,
or $e_1'$ contains a white vertex in $S^2-(D_1\cup D_2)$.

If $e_1'=e_2'$,
then there exists a lens.
This contradicts Lemma~\ref{NoLens}.
If $e_1'=e_2''$,
then this contradicts Lemma~\ref{LemmaTypeGH}(a).
If $e_1'=e_3'$,  
then this contradicts Lemma~\ref{TypeGE1E1}.
If $e_1'\ni w_4$,
then this contradicts Lemma~\ref{TypeGE1W4}.
Therefore, the edge $e_1'$ contains 
a white vertex in $S^2-(D_1\cup D_2)$.
Hence Claim~$1$ holds.
{\hfill {$\square$}\vspace{1.5em}}

{\bf Claim~2.}
The edge $e_1''$ contains a white vertex in $S^2-(D_1\cup D_2)$.

{\it Proof of Claim~$2$.}
By Condition~(e) of this section,
the edge $e_1''$ is not a terminal edge.
Thus there are five cases:
$e_1''=e_2'$, $e_1''=e_2''$, $e_1''=e_3'$, $e_1''\ni w_4$,
or $e_1''$ contains a white vertex in $S^2-(D_1\cup D_2)$.

If $e_1''=e_2'$,
then this contradicts Lemma~\ref{LemmaTypeGH}(a).
If $e_1''=e_2''$,
then there exists a lens.
This contradicts Lemma~\ref{NoLens}.
If $e_1''=e_3'$,
then this contradicts  Lemma~\ref{LemmaTypeGH}(b).
If $e_1''\ni w_4$,
then this contradicts  Lemma~\ref{TypeGE1W4}.
Therefore, the edge $e_1''$ contains 
a white vertex in $S^2-(D_1\cup D_2)$.
Hence Claim~$2$ holds.
{\hfill {$\square$}\vspace{1.5em}}

By Lemma~\ref{TypeGE5E5},
both of $e_5',e_5''$
contain white vertices in $S^2-(D_1\cup D_2)$.
Thus by Claim~1 and Claim~2,
each of the four edges
$e_1',e_1'',e_5',e_5''$
contains a white vertex in $S^2-(D_1\cup D_2)$.
Let $v_1,v_2,v_3,v_4$ be white vertices in
$S^2-(D_1\cup D_2)$ with
$v_1\in e_1',v_2\in e_1'',v_3\in e_5',v_4\in e_5''$.
Then the condition $w(\Gamma)=7$ implies
that 
\begin{enumerate}
\item[(1)] $w(\Gamma\cap(S^2-(G\cup D_1\cup D_2)))=2$.
\end{enumerate}
Hence, the set $\{v_1,v_2,v_3,v_4\}$ consists
of two white vertices.

Now, the four edges
$e_1',e_1'',e_5',e_5''$ are internal edges of label $m+1$
oriented inward at $w_1,w_1,w_5,w_5$,
respectively.
Thus, $e_1',e_1'',e_5',e_5''$ are oriented
outward at $v_1,v_2,v_3,v_4$,
respectively.
Moreover,
the five edges $e_2',e_2'',e_3',e_4',e_4''$
are oriented outward at $w_2,w_2,w_3,w_4,w_4$,
respectively.
Furthermore, we can show that
none of the nine edges 
$e_1',e_1'',e_5',e_5'',e_2',e_2'',e_3',e_4',e_4''$ 
are terminal edges.
Hence by IO-Calculation with respect to $\Gamma_{m+1}$
in $Cl(S^2-(D_1\cup D_2))$,
we have $w(\Gamma\cap(S^2-(D_1\cup D_2)))\geqq3$
(see Fig.~\ref{fig36}(b)).
This contradicts (1).
Therefore,
the graph $\Gamma_m$ does not contain
the graph as shown in Fig~\ref{fig10}(g).
\end{Proof}

\begin{lemma}
\label{MainTheoremChartAppX} 
{\rm (\cite[Theorem 1.1]{ChartAppX})}
Let $\Gamma$ be a minimal chart of type $(m;5,2)$.
Suppose that there exists a connected component of
$\Gamma_m$ with exactly five white vertices.
Then $\Gamma_m$ contains 
 one of the two graphs as shown in 
Fig.~\ref{fig10}$($g$)$,$($h$)$.
\end{lemma}

Now, we shall show the main theorem.

{\it Proof of Theorem~\ref{MainTheorem}.}
Suppose that there exists a minimal chart of $\Gamma$
of type $(m;5,2)$.

Suppose that there exists a connected component $G$
of $\Gamma_m$ with $w(G)=5$.
Then by Lemma~\ref{MainTheoremChartAppX},
the graph $\Gamma_m$ contains one of the two graphs as shown in 
Fig.~\ref{fig10}(g),(h).
However, this contradicts
Proposition~\ref{NoTypeH} and
Proposition~\ref{NoTypeG}.
Thus there exist at least two connected components
$G_1,G_2$ of $\Gamma_m$
with $w(G_1)\geqq1$ and $w(G_2)\geqq1$.

Now, by Lemma~\ref{LemmaNoLoop},
the chart $\Gamma$ does not contain any loop.
Hence $\Gamma_m$ does not contain any loop.
Thus by Lemma~\ref{GammaMFiveWhite}(b),(c),
the graph $\Gamma_m$ contains
a $\theta$-curve or an oval.
However, this contradicts
Lemma~\ref{NoThetaCurve} and Proposition~\ref{NoOvalSkewTheta}.
Therefore, there does not exist a minimal chart of $\Gamma$
of type $(m;5,2)$.
{\hfill {$\square$}\vspace{1em}}




\vspace{5mm}

\begin{minipage}{65mm}
{Teruo NAGASE
\\
{\small Tokai University \\
4-1-1 Kitakaname, Hiratuka \\
Kanagawa, 259-1292 Japan\\
\\
nagase@keyaki.cc.u-tokai.ac.jp
}}
\end{minipage}
\begin{minipage}{65mm}
{Akiko SHIMA 
\\
{\small Department of Mathematics, 
\\
Tokai University
\\
4-1-1 Kitakaname, Hiratuka \\
Kanagawa, 259-1292 Japan\\
shima@keyaki.cc.u-tokai.ac.jp\\
shima-a@tokai.ac.jp
}}
\end{minipage}

\newpage

\vspace{0.7cm}

{\bf List of terminologies}\vspace{2mm}\\
{\small $
\begin{array}{ll||}
\text{$k$-angled disk} & p5 \\
\text{boundary arc pair $(\alpha,\beta)$} & p11 \\
\text{BW-vertex} & p7 \\
\text{C-move equivalent} & p3 \\
\text{chart} & p3 \\
\text{complexity $(w(\Gamma),-f(\Gamma))$} & p4 \\
\text{feeler} & p5 \\
\text{free edge} & p3 \\
\text{hoop} & p4 \\
\text{internal edge} & p9 \\
\text{inward} & p3 \\
\text{inward arc} & p16 \\
\text{IO-Calculation} & p16 \\
\text{keeping $X$ fixed} & p6 \\
\text{lens} & p17 \\
\text{locally minimal} & p6 \\
\text{loop} & p8 \\
\text{M4-pseudo chart} & p14 \\
\end{array}
~~
\begin{array}{ll}
\text{middle arc} & p3 \\
\text{middle at $v$} & p3 \\
\text{minimal chart} & p4 \\
\text{outward} & p3 \\
\text{outward arc} & p16 \\
\text{oval} & p 8\\
\text{point at infinity $\infty$} & p4 \\
\text{pseudo chart} & p5 \\
\text{ring} & p4 \\
\text{RO-family} & p5 \\
\text{simple hoop} & p4 \\
\text{skew $\theta$-curve} & p8 \\
\text{special $k$-angled disk} & p5 \\
\text{terminal edge} & p4 \\
\text{type $(m;n_1,n_2,\cdots,n_k)$ for a chart} & p2 \\
\text{$\theta$-curve}& p8\\
\text{$(D,\alpha)$-arc} & p11 \\
\text{$(D,\alpha)$-arc free} & p11\\
\end{array}
$}

\vspace{0.5cm}

{\bf List of notations}\vspace{2mm}\\
{\small $
\begin{array}{ll||}
\text{$\Gamma_m$} & p2 \\
\text{$w(\Gamma)$} & p4 \\
\text{$f(\Gamma)$} & p4 \\
\text{${\rm Int}X$} & p5 \\
\text{$\partial X$} & p5 \\
\text{$Cl(X)$} & p5 \\
\end{array}
$~~
$\begin{array}{ll}
\text{$\partial \alpha$} & p5 \\
\text{${\rm Int}\alpha$} & p5\\
\text{$w(X)$} & p5 \\
\text{$c(X)$} & p6\\
\text{$a_{ij},b_{ij}$} & p14 \\
& \\
\end{array}
$
}

\end{document}